\documentclass{article} 
\usepackage{amsfonts}
\usepackage{amsopn}
\usepackage{epsfig}
\usepackage{amsmath}
\begin{document}

\newtheorem{theorem}{Theorem}[section]
\newtheorem{proposition}[theorem]{Proposition}
\newtheorem{lemma}[theorem]{Lemma}
\newtheorem{definition}[theorem]{Definition}
\newtheorem{remark}{Remark}[section]
\newtheorem{corollary}[theorem]{Corollary}

\renewcommand{\P}{\mathbb{P}}
\newcommand{\C}{\mathbb{C}}
\renewcommand{\O}{\mathcal{O}}
\newcommand{\G}{\mathbb{G}}
\renewcommand{\a}{\sum_{i=1}^a \alpha_i}
\renewcommand{\b}{\sum_{j=1}^b \beta_j} 
\newcommand{\X}{\mathcal{X}}
\newcommand{\Y}{\mathcal{Y}}
\newcommand{\Z}{\mathbb{Z}}
\newcommand{\M}{\overline{M}_{0,n}}
\newcommand{\F}{\mathbb{F}}

\title{Degenerations of Del Pezzo Surfaces and Gromov-Witten
  Invariants of the Hilbert Scheme of Conics}
\author{Izzet Coskun} \date{January 13, 2004}
\maketitle
\noindent {\bf Abstract:} This paper investigates low-codimension
degenerations of Del Pezzo surfaces. As an application we determine
certain characteristic numbers of Del Pezzo surfaces. Finally, we
analyze the relation between the enumerative geometry of Del Pezzo
surfaces and the Gromov-Witten invariants of the Hilbert scheme of
conics in $\P^N$.

\tableofcontents
\bigskip
\bigskip

\noindent 2000 {\it Mathematics Subject Classification}:  
14N15, 14N25, 14N35

\newpage
\setcounter{page}{1}

\section{Introduction:}

This paper investigates the degenerations of Del Pezzo surfaces $D_n$
embedded in $\P^N$ by their anti-canonical bundle. Due to the vast
number of possibilities, we restrict our attention to describing
simple specializations of $D_n$. As an application we determine some
characteristic numbers of Del Pezzo surfaces. Finally, we discuss the
relation between these numbers and the Gromov-Witten invariants of the
Hilbert scheme of conics. We work exclusively over the complex number
field $\C$.

This is a sequel to \cite{coskun1:degenerations} where we studied the
enumerative geometry of rational normal surface scrolls. Already in
that case, to obtain recursive formulae for the number of surfaces
incident to general linear spaces, we needed to impose strong
non-degeneracy assumptions by requiring enough of the linear spaces to
be points.  The case of Del Pezzo surfaces is more complicated, but
instructive to consider.  The new features of this case can be
summarized as follows:

\begin{enumerate}
\item Reducible surfaces that are limits of one-parameter families of
  scrolls are again unions of scrolls.  Del Pezzo surfaces exhibit a
  much larger variety of degenerations (\S 3). For example, a Del
  Pezzo surface can degenerate to a union of scrolls, a union of a
  Veronese surface and a scroll, a union of a Del Pezzo surface of
  lower degree and planes, a union of a rational cone and an elliptic
  cone. This partial list indicates that we cannot hope for a
  reasonable recursive formula for characteristic numbers of Del Pezzo
  surfaces via degeneration methods except in very special cases.

\item The hyperplane sections of Del Pezzo surfaces are not rational,
  but elliptic curves. The case of genus one curves in $\P^N$ is the
  last case where we have a firm understanding of the enumerative
  geometry of curves satisfying incidences with linear spaces
  (\cite{vakil:rationalelliptic}).  Consequently, the Del Pezzo
  surfaces lie at the perimeter of surfaces whose enumerative geometry
  we can analyze via degenerations.
  
\item Unlike scrolls Del Pezzo surfaces can have non-trivial moduli.
  This often makes it challenging to recognize limits.
\end{enumerate}

An interesting observation resulting from our investigations is that
degenerations of higher dimensional varieties exhibit qualitative
behavior fundamentally different from that of curves. Degenerations of
incidence and tangency conditions on curves with respect to linear
spaces result in a closed system of enumerative problems. The limits
of curves again satisfy similar conditions. The limits of surfaces, on
the other hand, can be subject to arbitrarily complicated conditions.
Since degeneration arguments are prevalent in algebraic geometry, this
crucial difference is important to note. \smallskip

We now sketch an outline of the paper. 

\noindent {\bf Notation.} Let $S^*$ be the dual of the tautological
bundle over $\G(2,N)$, the Grassmannian of planes in $\P^N$. Let $X^N$
denote $\P (\mbox{Sym}^2 S^*)$. $X^N$ is a projective bundle over
$\G(2,N)$. We can interpret it as the space of pairs of a plane and a
conic in the plane. \smallskip

\noindent {\bf The limits.} We produce a list of potential
non-degenerate limits of $D_n$ that can occur in one-parameter
families (\S 3) using a classical theorem of Del Pezzo and Nagata (\S
2.3), which classifies surfaces of degree $n$ in $\P^n$.  We then exhibit
families realizing the degenerations of $D_n$ relevant to our counting
problems and we describe the limiting positions of geometrically
significant curves.

We use two techniques to construct families of $D_n$ specializing to a
given limit. We specialize the base points of the linear system of
cubics on $\P^2$ in various ways to obtain classical constructions.
More interestingly, since Del Pezzo surfaces $D_n$ are ruled by
conics, we can interpret them as curves in $X^N$. Let $d_n$ denote the
cohomology class of a curve in $X^N$ arising from a one-parameter
family of conics on $D_n$. Given a potential limit surface, we can try
to find a curve $C$ of conics in the class $d_n$ which sweeps it. If
we can deform $C$ to a curve of conics arising from a smooth $D_n$,
then we can conclude that the surface arises as a degeneration of
$D_n$. \smallskip

\noindent {\bf Example.} For instance, it takes ingenuity to find a
specialization of the base points in order to obtain a family of $D_n$
$(n<8)$ degenerating to the projection of the rational scroll
$S_{2,n-2}$ from a point on the plane of a conic on $S_{2,n-2}$.
However, it is easy to see that the curve of reducible conics
consisting of a fiber line and the double line deforms to a curve of
conics on $D_n$ (see \S 3). \smallskip

\noindent {\bf Characteristic numbers.} By the {\it characteristic
  number problem} we mean the problem of computing the number of
varieties of a given type that meet the `appropriate' number of linear
spaces in general position. Classically characteristic numbers also
allow tangency conditions; however, in this paper we will consider
only incidence conditions. Using our description of the degenerations
of $D_n$ we determine some characteristic numbers of $D_3$ and $D_4$.
Although most of these numbers can also be obtained by classical
methods, our method has the advantage of circumventing tedious
cohomology calculations and yields numbers of surfaces satisfying
divisorial conditions, which are hard to obtain classically. We can
also determine a few of the characteristic numbers of $D_5$.  However,
for $n \geq 6$ and essentially for $n=5$, the degenerations get too
complicated for the method to terminate and give actual numbers. If
instead we ask for the number of $D_n$ containing a fixed degree $n$
elliptic normal curve and satisfying incidences with linear spaces,
then the degeneration method gives a few more answers for $n=5$ and 6.
I do not know of a classical method to compute these numbers when
$n>4$.  \smallskip

\noindent {\bf The method.} In order to count surfaces incident to
various linear spaces, we degenerate the linear spaces one by one to a
hyperplane $H$ until we force any surface meeting them to become
reducible. If we have enough point conditions to satisfy our
non-degeneracy assumptions, we know the possible reducible surfaces
that occur in the limit. We can hope to count surfaces by further
breaking each of the pieces of the limit surface to obtain simpler
surfaces. This hope is in general upset by the appearance of singular
surfaces and more and more complicated conditions on hyperplane
sections of the surfaces. However, the method still works in many
cases (see examples in \S 4) and with some effort should extend to
more cases than covered here.

\smallskip

\noindent {\bf Gromov-Witten invariants of $X^N$.} It is natural to
ask for the relation between the enumerative numbers for $D_n$ and the
Gromov-Witten invariants of $X^N$. $X^N$ is not a convex variety (\S
2.2), i.e.\  there are maps $f: \P^1 \rightarrow X^N$ for which
$h^1(\P^1, f^* T_{X^N}) \not= 0$. Its Kontsevich spaces of genus zero
stable maps often have components of more than the expected dimension.

When a variety $V$ is homogeneous (in particular convex), then the
Gromov-Witten invariants count the number of curves that meet general
subvarieties of $V$. However, when the variety is not convex, there
can be virtual contributions to the Gromov-Witten invariants. It is
usually a hard problem to decide when the Gromov-Witten invariants of a
non-convex space are enumerative. Gathmann \cite{gathmann:blowup} and
G\"ottsche and Pandharipande \cite{gottschepan:blowup} discuss this
problem for the blow-ups of $\P^N$ and $\P^2$, respectively.

In general the Gromov-Witten invariants of $X^N$ for the class $d_n$
are not enumerative (\S 6). However, we prove that the Gromov-Witten
invariants involving incidences to linear spaces are enumerative when
$n=3$ and when $n=4$ provided that there are not any $\P^3$s incident
to all the linear spaces. As a corollary we compute some Gromov-Witten
invariants of $X^N$.  \smallskip

\noindent {\bf Acknowledgements.} It is a pleasure to thank Brendan
Hassett, Ciro Ciliberto, Mihnea Popa, Jason Starr, Dan Avritzer and
Maryam Mirzakhani for fruitful discussions. I am especially grateful
to Ravi Vakil and my advisor Joe Harris for their ideas and invaluable
suggestions during the course of this project. I would also like to
thank Jun Li, the Stanford Mathmatics Department and especially Ravi
Vakil for their hospitality.

\section{Preliminaries}

\subsection{ Del Pezzo surfaces}
In this subsection we discuss the basic geometry of Del Pezzo
surfaces.  For more details consult \cite{beauville:surface} Ch. 4,
\cite{joe:thebook} \S 1 Ch. 4  or \cite{friedman:surfaces} Ch. 5.  \smallskip

\noindent  {\bf Del Pezzo surfaces} are smooth complex surfaces with
ample anti-canonical bundle $-K$. Except for $\P^1\times\P^1$, they
can be realized as the blow-up of $\P^2$ in fewer than 9 points no
three of which lie on a line and no six of which lie on a conic. To
have a more uniform discussion we exclude $\P^1 \times \P^1$.  We
denote Del Pezzo surfaces by $D_n$ where $n$ is the degree $K^2$ of
the anti-canonical bundle.  Equivalently, $D_n$ is the blow-up of
$\P^2$ in $9-n$ general points $p_i$.  The anti-canonical series
$\left| -K \right|$ on a Del Pezzo surface can be interpreted as the
linear series of cubics on $\P^2$ having $p_i$ as base points,
therefore

\begin{displaymath}
h^0(D_n, -K) = 10- n.
\end{displaymath}
We limit our discussion to Del Pezzo surfaces $D_n$ embedded in $\P^n$
by their anti-canonical
bundle, i.e.\  to $D_n$ with $n \geq 3$. \smallskip

\noindent {\bf Geometric description.} These surfaces display a rich
geometry and often have nice determinantal descriptions.  $D_3$ is a
cubic surface in $\P^3$. $D_4$ is the complete intersection of two
quadric threefolds in $\P^4$.  $D_5$ is a fourfold hyperplane section
of the Grassmannian $\G(1,4)$ under its Pl\"ucker embedding.  $D_6$ is
a two fold hyperplane section of the Segre embedding of $\P^2 \times
\P^2$ or it is the hyperplane section of the Segre embedding of $\P^1
\times \P^1 \times \P^1$. Finally, $D_9$ is the cubic Veronese
embedding of $\P^2$ in $\P^9$.  \medskip

\noindent {\bf The Picard group of $D_n$} is isomorphic to $\Z^{10-n}$
generated by the classes $H$, the pull back of the hyperplane class
from $\P^2$, and $E_i$, $1 \leq i \leq 9-n$, the exceptional divisors
of the blow-up. The intersection pairing is

\begin{displaymath}
H^2 = 1,  \ \ \ \ \ \ H \cdot E_i = 0, \ \ \ \ \ \ E_i \cdot E_j = -
\delta_{i,j} .
\end{displaymath}  
In terms of these classes the anti-canonical class is
$-K = 3H - \sum_{i=1}^{9-n} E_i.$

By Bertini's theorem a general hyperplane section of $D_n$ is a smooth
elliptic curve of degree $n$. These curves are projectively normal.
\smallskip

\noindent {\bf Curves on $D_n$.} During the degenerations it is
important to know the limits of lines, conics and hyperplane sections
on $D_n$. On $D_n$ the effective curve classes containing an
irreducible curve of a given arithmetic genus $g$ and degree $d$ are
easy to determine. We can express the class of any curve as $aH
-\sum_{i=1}^{9-n} b_i E_i$.  Since the surface is embedded by $\left|
  -K \right|$, the degree condition implies that
$$3a-\sum_{i=1}^{9-n}b_i = d.$$
The genus formula translates to $$a^2
- \sum_{i=1}^{9-n} b_i^2=2g-2+d.$$
Using the Cauchy-Schwarz inequlity
$$\left( \sum_{i=1}^{9-n} b_i \right)^2 \leq (9-n)\sum_{i=1}^{9-n}
b_i^2,$$
we find the choices for $a$ and then solve for the $b_i$
satisfying the two equations. In the rest of the paper we will use
this scheme to determine curve classes without further mention. For
the convenience of the reader we enumerate the classes of lines and
conics on $D_n$.

\begin{lemma}\label{lines}
  On the Del Pezzo surfaces $D_n$ ($n\geq 3$) the classes of lines are
  $E_i$, $1 \leq i \leq 9-n$, $H-E_i-E_j$, $i \not= j$ and $2H - E_a-
  E_b-E_c-E_d - E_e$ where $a,b,c,d,e$ are distinct, whenever these
  classes exist.
\end{lemma}

\begin{lemma}\label{conics}
  On the Del Pezzo surfaces $D_n$ ($n \geq 3$) the classes of conics
  are $H - E_i$, $2H - E_a - E_b-E_c- E_d$, $3H - 2E_a
  -E_b-E_c-E_d-E_e-E_f$ where $ a, b, c, d, e, f$ are distinct,
  whenever these classes exist.
\end{lemma}

\noindent {\bf Moduli of Del Pezzo surfaces.} The surfaces $D_3$ and
$D_4$ have a four and two dimensional moduli space, respectively. We
will not be concerned with the construction or properties of these
moduli spaces.   \smallskip

\noindent {\bf Singular Del Pezzo surfaces.} A {\it singular Del Pezzo
  surface} $D_n^{(s)}$ is an irreducible surface of degree $n$ in
$\P^n$ which has isolated double points.  $D_n^{(s)}$ is also the
image of the blow-up of $\P^2$ in $9-n$ points. $D_n^{(s)}$ arises
when the points we blow up to obtain $D_n$ become infinitely near,
fail to be in general linear position or lie on a conic (when $n=3$).

The list of the combinations of double points that occur on
$D_n^{(s)}$ is long.  When $n=3$, Bruce and Wall give a very nice
description \cite{bruce:wall:cubic}. The type and combination of the
double points that occur on a cubic surface are all obtained by
deleting vertices (and the edges adjacent to them) from the extended
$E_6$ diagram.

\begin{figure}[htbp]
\begin{center}
\epsfig{figure=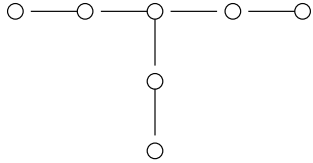}
\end{center}
\caption{The extended $E_6$ diagram}
\label{Figure 1}
\end{figure}
\noindent Conversely, every combination of Du Val singularities that arises by
deleting vertices from the extended $E_6$ diagram occurs on some cubic
surface.  \smallskip

\noindent {\bf Tangent planes to $D_n$ along curves.} When studying
limiting positions of hyperplane sections in one-parameter families of
surfaces, it is essential to have estimates on the dimension of the
space of hyperplanes tangent to the components of the limit surface
along common curves. The dimension of the space of hyperplanes tangent
to a smooth $D_n \subset \P^n$ along a line is $\max(-1, n-5)$.
However, if the surface has a double point this estimate can change.
For example, on a smooth $D_4 \subset \P^4$ there are not any
hyperplanes tangent to the surface everywhere along a line; however,
when the surface acquires an ordinary double point, there can be one.
Consequently, one has to exercise caution to consider all possible
singularities when giving estimates. For simplicity, we will often
exclude singular surfaces from the discussion. \smallskip

\noindent {\bf The dimension of the space of Del Pezzo surfaces in
  $\P^N$.} If we fix linear spaces $\Lambda^{a_i}$ of dimension $a_i$
in general position such that $$\sum_i \left( N-2-a_i \right) = N(n+1)
-n + 10,$$
then there will be finitely many smooth $D_n$ meeting all
the linear spaces by the following dimension count.

\begin{lemma}\label{dimdelpezzo}
  The dimension of the locus in the Hilbert scheme whose general point
  corresponds to a smooth $D_n$ in $\P^N$
  is
\begin{displaymath}
N(n+1) -n + 10.
\end{displaymath}
\end{lemma}

\noindent {\bf Proof:} Realize $D_n$ as the image of a map from the
blow-up of $\P^2$ at $9-n$ points by choosing $N+1$ sections in
$\left| -K \right|$
and projectivizing. We need to add the dimension of the moduli space
or subtract the dimension of the automorphism group which amounts to
adding $10-2n$.  $\Box$ \smallskip

We will determine the number of $D_n$ in some cases using
degenerations. For future reference we recall the following well-known
fact (see \cite{vakil:rationalelliptic} \S 5).
\begin{lemma}\label{dimellipticcurve}
  The dimension of the component of the Hilbert scheme whose general
  point corresponds to a smooth elliptic curve of degree $n+1$
  spanning a $\P^n$ in $\P^N$ is $$(N-n)(n+1) + (n+1)^2.$$
\end{lemma}

\noindent {\bf Rational Scrolls.} During the degenerations of $D_n$
we will encounter rational surface scrolls. We refer the reader to
$\cite{coskun1:degenerations}$ for a detailed discussion of their
geometry.

A rational normal scroll $S_{k,l}$ is abstractly the Hirzebruch
surface $F_{l-k}$ embedded in $\P^{k+l+1}$ by the complete linear
series $e+lf$, where $e,f$ are the usual generators of the Picard group of
$F_{l-k}$ satisfying $e^2=k-l, \ f^2=0, \ e \cdot f=1$. The classes
$e,f$ are the classes of the exceptional curve $E$ and of a fiber $F$,
respectively. The surface can be explicitly constructed by taking a
rational normal $k$ curve and a rational normal $l$ curve with
disjoint linear spans; choosing an isomorphism between the curves; and
taking the union of lines joining corresponding points.

We will refer to a curve class $e+mf$ as a {\bf section class} and to
a curve class $2e + mf$ as a {\bf bisection class}. Irreducible curves
in section and bisection classes are sections and bisections of the
projective bundle over $\P^1$, respectively. In addition to the
cohomology calculations in \S 2 of \cite{coskun1:degenerations}, we
will use

\begin{displaymath}
  h^0
(F_r, \O_{F_r} (2e+(r+2)f)) = \left\{ \begin{array} {c@{\quad:\quad}l}
    9 & r \leq 2 \\ h^0 (F_r, \O_{F_r} (e+(r+2)f)) & r 
    \geq 3 \end{array} \right. 
\end{displaymath}
which follows by considering the exact sequences
\begin{equation*}
\begin{split}
  & 0 \rightarrow \O_{F_r} (2e+mf) \rightarrow \O_{F_r} (2e+(m+1)f)
  \rightarrow \O_{F} (2) \rightarrow 0 \\
  & 0 \rightarrow \O_{F_r} (e+mf) \rightarrow \O_{F_r} (2e+mf) \ 
  \rightarrow \O_{E} (m-2r) \rightarrow 0.
\end{split}
\end{equation*}

\noindent {\bf The Veronese surface.} During the degenerations
we will also encounter the Veronese surface, the image of $\P^2$ in
$\P^5$ given by the complete linear system of conics. A central fact
is that the Veronese surface together with the rational normal scrolls
are the only non-degenerate irreducible surfaces of degree $n-1$ in
$\P^n$ (\cite{joe:thebook} p. 525).

\subsection{The Geometry of the Space of Conics}

Since we will rely on the description of $D_n$ as a curve of conics,
we recall the basic facts about the Hilbert scheme of conics in
$\P^N$. \smallskip

Let $S^*$ denote the dual of the tautological bundle on the
Grassmannian $\G(2,N)$ of planes in $\P^N$. Recall that $X^N$
was defined to be $\P(\mbox{Sym}^2 S^*)$. Let
$$\pi:X^N \rightarrow \G(2,N)$$
be the natural projection map. 

The Chow ring of $X^N$ is generated by the pull-back of the classes on
$\G(2,N)$ and the first chern class of the tautological bundle on
$X^N$ (\cite{fulton:intersection} 8.3.4). The Picard group has 2
generators which we can express in terms of geometric cycles. Let
$\omega$ be the class of conics which meet a fixed $\P^{N-2}$ (i.e.\ 
the pull-back of the hyperplane class of $\P^N$ by the natural
morphism) and $\eta$ the class of conics whose planes meet a fixed
$\P^{N-3}$ (more precisely, the first chern class of the tautological
bundle on $X^N$). The locus of reducible conics $\Delta$, whose class
we denote by $\delta$, is also a divisor.  The intersection products
below (see Table 1) imply that $\delta = 3 \omega - 4 \eta$.

Similarly, there are two natural curve classes on $X^N$. Let $a$ be
the class of conics that lie in a fixed plane and pass through 4
general points on the plane. Let $b$ be the class of conics that are
cut out on a fixed quadric surface by a pencil of hyperplanes. Each
conic class contains a one-parameter family of disjoint conics on
$D_n$, $n\not= 9$, so gives rise to a curve in $X^N$. For example,
$D_3$ gives rise to 27 rational curves in $X^3$.  Let $d_n$ denote the
class in $X^N$ of the curve of conics on $D_n$ in a fixed cohomology
class. Table 1 below implies that $$d_n = (4-n) a + (n-2) b.$$
\smallskip

\centerline{ \begin{tabular}{l|c|c|c}
$\cdot$ & $\omega$ & $\eta$ & $\delta$ \\ \hline
$a$ & 1 & 0 & 3 \\ \hline
$b$ & 2 & 1 & 2 \\ \hline
$d_n$ & $n$ & $n-2$ & $8-n$
\end{tabular}}
\smallskip

\centerline{Table 1: Intersection Products}
\medskip

For $\beta \in H^2(X^N, \Z)$ let $\overline{M}_{0,m}(X^N, \beta)$
denote the Kontsevich space of stable maps in the class $\beta$. In
general the Kontsevich spaces of $X^N$ can have components of larger
than expected dimension. One of the simplest examples is
$\overline{M}_{0,0}(X^N, -d a + db)$ for $d >1$.  Using the Euler
sequence for the tangent bundle (\cite{fulton:intersection} 3.2.11)
one can check that the canonical bundle of $X^N$ is given by $$K_{X^N}
= -6 \omega + (7-N) \eta.$$
Hence, the expected dimension of
$\overline{M}_{0,0}(X^N, -d a + db)$ is
$$\dim\ X^N + c_1(X)(-d a + db) -3 = Nd-d+ 3N-4.$$
Take a cone $C$
over a rational curve of degree $d$ and a line $l$ not necessarily
contained in $C$, but meeting it at the vertex. The curve in $X^N$
whose points correspond to the union of $l$ with a line of $C$ has
class $-da + db$. The dimension of cone and line pairs is $Nd+3N-5$.
So when $d >1$, this provides us with a Kontsevich space of the
``wrong'' dimension. We will see more examples in \S 6. For future
reference we note that the expected dimension of $\overline{M}_{0,m}
(X^N,d_n)$ is
$$N(n+1) - n + 10 + m.$$
When $m=0$, this agrees with the dimension in Lemma \ref{dimdelpezzo}.

\subsection{The Classification of degree $n$ surfaces in $\P^n$.}
In this subsection we state the classification theorem for reduced,
irreducible, non-degenerate surfaces of degree $n$ in $\P^n$. This is
a classical theorem of Del Pezzo and Nagata whose proof can be found
in \cite{nagata:memoir}.

\begin{theorem}\label{classification}
  An irreducible, reduced, non-degenerate surface of degree $n$ in
  $\P^n$ is one of the following:
  
1. A projection to $\P^n$ of a scroll of degree $n$ in $\P^{n+1}$,

2. A projection to $\P^4$ of the Veronese surface in $\P^5$,

3. A Del Pezzo surface, possibly with finitely many isolated double
points,

4. The image of $F_0$ or $F_2$ in $\P^8$ given by their anti-canonical map,

5. A cone over an elliptic curve of degree $n$ in $\P^{n-1}$. 
\end{theorem}

\section{The limits of Del Pezzo surfaces}

In this section we describe the non-degenerate surfaces in $\P^n$ that
can arise as limits of $D_n$. We appeal to the description of $D_n$ as
a curve in $X^n$. Since $D_9$ does not contain any conics, we restrict
the values of $n$ to $3 \leq n \leq 8$.  \smallskip

\subsection{Constraints on the Degenerations} 

\noindent {\bf Notation.} Let $f: \Y \rightarrow B$ be a flat family
of surfaces in $\P^N$ over a smooth, connected curve. Let $b_0 \in B$
be a marked point and let $\Y_0$ denote the fiber of $f$ over $b_0$.
We assume that $\Y_b$ for points $b \not= b_0$ is a Del Pezzo surface
$D_n$. We also assume that $\Y_0$ spans a $\P^n$ and its components
are reduced. We preserve the notation we used in \S 2 for $X^N$.
\smallskip 

\noindent {\bf Adjacent components.} Let $Y$ be a reducible surface
connected in codimension 1. We refer to two components that share a
common curve as {\it adjacent components.} The {\it dual graph} of the
surface consists of a vertex for each irreducible component and an
edge between adjacent ones.

\begin{lemma}\label{tree}
  $\Y_0$ is a surface of degree $n$ whose components are ruled by
  lines or conics. There exists a connected subgraph containing all
  the vertices of its dual graph such that adjacent components share a
  common line or conic.
\end{lemma}

\noindent {\bf Proof:} Each $D_n$ in the family gives rise to a
collection of rational curves in $X^N$. After a finite base change
totally ramified over $b_0$, we can select a conic class on each
surface away from $b_0$. We denote the new family by $\Y ' \rightarrow
B'$. This family induces a curve in the Kontsevich space of stable
maps $\overline{M}_{0,0}(X^N,d_n)$.  The limit of the family in
$\overline{M}_{0,0}(X^N,d_n)$ is a map from a tree of rational curves
to $X^N$.  The restriction of the universal curve over $X^N$ to the
family of curves maps to $\P^N$ giving rise to a family of surfaces
which agrees with $\Y'$ except possibly over $b_0'$. There is a scheme
structure on the limit surface which makes the family flat. Since over
a smooth curve there is a unique way to complete a family to a flat
family, this family agrees with our original family. We conclude that
the components of $\Y_0$ are ruled by conics or lines. The last
assertion is clear. $\Box$

\begin{proposition} \label{junk}
  Each component of $\Y_0$ is one of the following: 
  
  1. A Veronese surface in $\P^5$ or a non-degenerate scroll of degree
  $k$ in $\P^{k+1}$,
  
  2. A projection of one of the surfaces in 1 to a surface with a
  double line in a one dimensional lower projective space,
  
  3. A cone over an elliptic curve of degree $k$ in $\P^{k-1}$, or
  
  4. A Del Pezzo surface $D_n$, possibly with isolated double points.
\end{proposition}

\noindent {\bf Proof:} If $\Y_0$ is irreducible, then it is a
non-degenerate surface of degree $n$ in $\P^n$. $\Y_0$ cannot be the
anti-canonical image of $F_0$ or $F_2$ since these surfaces do not
contain any lines. The flat limit of the lines in the family of $D_8$
would be a line. Similarly, $\Y_0$ cannot be the projection of a
rational scroll or Veronese surface with isolated singularities. The
hyperplane section of such a surface has arithmetic genus 0 instead of
1. By Theorem \ref{classification} we conclude that $\Y_0$ is one of
the surfaces in cases 2, 3 or 4. \smallskip

\noindent {\bf Two components.} Suppose $\Y_0$ has two irreducible
components $W$ and $Z$ of degrees $d_W$ and $d_Z$.  By Lemma
\ref{tree} they share a line or a conic.  \smallskip

\noindent {\bf Suppose $W$ and $Z$ meet in a conic} (possibly reducible
or non-reduced).  Then the linear spaces they span contain a common
plane, so their total span is at most $\P^{d_Y+d_Z}$.  We conclude
that the surfaces must be minimal degree surfaces, so they are one of
the surfaces in case 1. Since the Veronese surface does not contain
any lines, at most one of the surfaces can be a Veronese surface.
Therefore, the surfaces are either two scrolls meeting along a conic
or a Veronese and a scroll meeting along a conic. \smallskip

\noindent {\bf Suppose $W$ and $Z$ meet along a line.} If both of the
surfaces are minimal degree surfaces and meet generically transversely
along the line, then their union cannot be a limit of Del Pezzo
surfaces since their hyperplane sections have arithmetic genus 0
instead of 1. It is possible for two rational cones tangent along a
line to be a limit of $D_n$, but this is included in the previous
case.

We can, therefore, assume that $W$ spans only a $\P^{d_W}$. $Z$
must span a $\P^{d_Z + 1}$ and meet $W$ along the line generically
transversely. Since a Veronese surface does not contain any lines, $Z$
is a rational normal scroll.  Theorem \ref{classification} implies
that $W$ is a Del Pezzo surface possibly with finitely many double
points, a cone over an elliptic curve or the projection of a scroll or
a Veronese.  Moreover, if $W$ is the projection of a scroll or the
Veronese, it must have a double line because otherwise a general
hyperplane section would have arithmetic genus 0.  \smallskip

\noindent {\bf Further constraints. } If $W$ is a Del Pezzo surface
and $Z$ is not a plane, then both $W$ and $Z$ must be singular.
Suppose to the contrary that $W$ is a smooth $D_k$ and $Z$ is a scroll
of degree greater than 1 meeting it in a line $l_C$.  The limit of a
curve of conics on $D_n$ is a connected curve of conics on $\Y_0$.
Since $W$ and $Z$ do not have a conic in common, $W$ is not ruled by
lines and the double of $l_C$ is not a conic class on $W$, the conics
on $\Y_0$ must consist of a curve of conics on $W$ union a curve of
lines on $Z$ together with a fixed line $l_F$ on $W$ intersecting each
line on $Z$. Hence, $Z$ must be a cone. The fixed line $l_F$ cannot be
$l_C$, so the conic class on $W$ is the conic class $[l_F] + [l_C]$.
The intersection of the variety of reducible conics with this curve
has more than $8-n$ isolated points (see table in \S 2), hence the
curve cannot be deformed to a smooth curve in the class $d_n$. We
conclude that $W$ is also singular. Furthermore, the same argument
shows that as the degree of the scroll increases, the residual Del
Pezzo surface is forced to have worse singularities.  For example, if
$D_n$ breaks into a cubic cone union $D_{n-3}$, then $D_{n-3}$ must
have a singularity worse than an ordinary double point. Using the list
of singularities and the number of lines on the singular surfaces it
is not too hard to make a list of possibilities. \smallskip

\noindent {\bf Remark.} In case $W$ is not $D_3$ or $D_4$ the previous
constraint follows by an elementary dimension count. However, since
$D_3$ and $D_4$ have moduli the dimension count only shows that $W$
cannot be a general smooth $D_3$ or $D_4$. \smallskip

By a similar argument if $W$ is a cone over an elliptic curve, then
$Z$ is a rational cone with matching vertex. The limit curve of conics
has to be reducible. One component of the curve must consist of line
pairs on $W$ joining the points identified by the hyperelliptic
involution on the elliptic curve. The other component must consist of
a fixed line in $W$ union lines in $Z$. The claim follows. \smallskip

\noindent {\bf More components.} Now we allow $\Y_0$ to have more than
two components. \smallskip

\noindent {\bf Observation.} If a subsurface $S$ of
$\Y_0$ of degree $d$ spans exactly $\P^d$, then all the remaining
components are minimal degree scrolls meeting the components adjacent
to them in lines.  The components adjacent to $S$ need to be attached
to $S$ along at least a line. To have the resulting surface be
non-degenerate, the surface we attach must have maximal possible span
for its degree and meet at most one of the components of $S$ in a
line. The observation follows by induction.

\noindent $\bullet$ Suppose $\Y_0$ contains three components
$U_i$ of degree $d_i$ pairwise meeting in distinct curves. $U_i$ spans
at most $\P^{d_i + 1}$. The three components together span at most a
linear space of dimension $\sum_i d_i$ with equality if and only if
each of the surfaces have maximal span and their common curves are
concurrent lines. By the observation we conclude that each component
is a scroll and the $U_i$ each contain a pair of intersecting lines.
\smallskip
        
\noindent $\bullet$ Suppose $\Y_0$ contains two components meeting in a
conic, possibly reducible or non-reduced. Then the two components are
either two scrolls meeting along a conic or a Veronese surface and a
scroll meeting along a conic. By the observation all the other
components must be scrolls.  \smallskip

\noindent $\bullet$ We can assume that all components of $\Y_0$ meet
pairwise in lines and no three components meet pairwise in distinct
curves. Suppose one of the components $U$ of degree $d$ spans $\P^d$.
By Theorem \ref{classification} and the argument given for the case
when $\Y_0$ has two components, $U$ is the projection of a Veronese
surface or a rational scroll with a double line, a cone over an
elliptic normal curve or a Del Pezzo surface possibly with finitely
many singularities. By the observation all the other components are
scrolls. By an argument similar to the two component case we can
deduce that if $V$ is a smooth Del Pezzo surface then all the adjacent
scrolls are planes and they are joined to $V$ along non-intersecting
lines. In case $V$ is a cone over an elliptic normal curve, the
adjacent components are rational cones whose vertices coincide with
the vertex of $V$.  \smallskip

\noindent $\bullet$ Finally, we can assume that all the components are
rational scrolls meeting pairwise in lines. The hyperplane sections of
the surface ought to have arithmetic genus 1. This concludes the
description of the components.  $\Box$

\begin{corollary}
1. At most one component of $\Y_0$ is a Veronese surface. If a
   component is a Veronese surface, all other components are rational
   normal scrolls.

2. If a component $U$ of $\Y_0$ of degree $d$ spans only $\P^d$, then all
   the other components are rational normal scrolls. If $U$ is a
   smooth Del Pezzo surface or a cone over an elliptic normal curve,
   the scrolls have to satisfy the constraints described in the proof
   of Proposition \ref{junk}.
\end{corollary}

\subsection{Explicit Families Realizing the Degenerations}

One might hope that few of the surfaces described in Proposition
\ref{junk} actually occur as components of the limtis of
$D_n$. Unfortunately this is not so. We now give examples of families
realizing most of the irreducible or two-component possibilities listed
in Proposition \ref{junk}. We assume $n<8$ throughout. 
\smallskip

\noindent {\bf Notation.} Let $A$ denote a  disk in $\C$. Let $a_0 \in
A$ denote a marked point in the disk. We denote sections of a family
of varieties over $A$ by $s_i$. Finally, let $N_{X/Y}$ denote the
normal bundle of $X$ in $Y$.
\smallskip

\noindent {\bf I.} The Del Pezzo surfaces with isolated double points
all arise by specializing the base points of the linear system of
cubics on $\P^2$.  From this description one can determine the
limits of the lines and conics.

For example, to construct a family of cubic surfaces specializing to a
cubic surface with four $A_1$ singularities, take $A \times \P^2$ and
6 general sections $s_i$ which specialize to the intersection points
of 4 lines $l_i$ in the $\P^2$ over $a_0$. Blow up the sections $s_i$
in $A \times \P^2$. The dual of the relative dualizing sheaf
restricted to each fiber embeds the fiber away form $a_0$ as a $D_3$
in $\P^3$.  The image of the fiber over $a_0$ is a cubic surface with
the required singularities. There are 9 lines on the limit surface:
the image of the 6 exceptional divisors and the image of the three
lines joining the pairs of points that do not lie on an $l_i$. The
descriptions of conics is similar.


\begin{lemma}\label{smoothing}
  Suppose $C$ is a smooth rational curve in $X^n$ contained in the
  locus of reducible conics $\Delta$. Suppose $C$ does not intersect
  the locus of non-reduced conics and $[C] \cdot \delta \geq 0$. Then
  $C$ can be deformed away from $\Delta$.
\end{lemma}

\noindent {\bf Proof:} $C$ is contained in the smooth locus of
$\Delta$. Away from the locus of non-reduced conics $\Delta$ is a
homogeneous variety and its tangent bundle is generated by global
sections. Consequently, $N_{C/\Delta}$ is generated by global
sections. Using the exact sequence
$$0 \rightarrow N_{C/\Delta} \rightarrow N_{C / X^n} \rightarrow \O_C
(\Delta) \rightarrow 0$$
and the assumption that $[C] \cdot \delta
\geq 0$, we conclude $h^1(C,N_{C/X^n}) = 0$ and that $H^0
(C,N_{C/\Delta})$ does not surject onto $H^0 (C,N_{C/X^n})$. Hence,
the first order deformations of $C$ are unobstructed and a general
first order deformation of $C$ does not lie in $\Delta$. The lemma
follows. $\Box$  \smallskip

\noindent {\bf II. Degeneration of $D_n$ to a scroll with a double
  line.}  The projection of a scroll $S_{1,2}$ or $S_{2,l}$, $l \leq
6$, from a general point on the plane of a conic in the surface has
the same Hilbert polynomial as a Del Pezzo surface. We will show that
these surfaces are limits of Del Pezzo surfaces. These scrolls can be
further degenerated to more unbalanced scrolls.

The union of the double line with the fibers gives rise to a curve $C$
of reducible conics in the Hilbert scheme $X^n$ contained in the
smooth locus of $\Delta$.  Since $\delta \cdot [C] = 8-n$, we conclude
by Lemma \ref{smoothing} that $C$ can be deformed away from $\Delta$.
The resulting curve has the same class as a curve arising from a Del
Pezzo surface.  The surface in $\P^n$ spanned by the conics is a
non-degenerate, irreducible surface of degree $n$ ruled by reduced and
generically irreducible conics. Since the dimension of scrolls with a
choice of curve of conics sweeping the surface once is smaller than
the dimension of the deformations of $C$, by Theorem
\ref{classification} we conclude that the surface is $D_n$ (recall
$n<8$). \smallskip

\noindent {\bf Alternative  construction when $n \leq 5$.}
For concreteness assume $n=5$. Take $A \times \P^2$ and specialize 4
general sections $s_i$ to the same point $p$ on the central fiber.
Blow up the total space at $p$, then along the proper transform of the
sections $s_i$.
Denote the total space of the resulting threefold by $X$. The central
fiber is the union of $F_1$ with a $\P^2$ blown up at 4 points. Denote
these two components of the central fiber by $F$ and $P$,
respectively. Take the linear system which restricts to $\left| e+3f
\right|$ on $F$ and to $2L - \sum_{i=1}^4 E_i$ on $P$, where $L$ is
the line class on $\P^2$ and $E_i$ are the exceptional divisors of the
blow-ups. The linear series on the $F_1$ component is not complete,
but must match the linear series on $P$. The latter contracts the
surface to a double line.  This constructs the desired degeneration of
$D_5$.  \smallskip

\noindent{\bf III. Degeneration of $D_n$ to a cone over an elliptic
  curve.} Every surface degenerates to a cone over a hyperplane
section (possibly with some embedded structure at the cone point) by
taking the limit of a one parameter family of projective
transformations fixing the hyperplane. Since the cone over an elliptic
curve of degree $n$ in $\P^{n-1}$ has the same Hilbert polynomial as
$D_n$ there are no embedded components in this case. By degenerating
the elliptic curve which is the base of the cone into an elliptic
curve with rational tails, one obtains degenerations of $D_n$ into an
elliptic cone union rational cones.  \medskip

\noindent {\bf IV. $D_4$ degenerates to the projection of a Veronese
  surface with a double line.} Both $D_4$ and the Veronese surface
with a double line are complete intersections of two quadric
threefolds. Since a general complete intersection is a $D_4$, to
obtain such a degeneration it suffices to specialize the quadrics. 

To see that a Veronese surface with a double line is a complete
intersection of quadric threefolds it suffices to observe that such a
Veronese surface is given by the map $$(x_0, x_1, x_2) \mapsto (x_0^2,
x_1^2, x_2^2, x_0x_1, x_0x_2) $$
in projective coordinates, hence two
quadric threefolds contain it.

A degeneration of $D_n$ into $D_4$ union other surfaces meeting $D_4$
along lines further degenerates to a surface where a component is a
Veronese with a double line.  \medskip

\noindent {\bf V. Degeneration of $D_n$ to $D_{n-1}$ union a plane.}
In $A \times \P^2$ blow up $9-n$ disjoint sections $p_i (a)$ of $A$ in
general position. In the fiber over $a_0$ blow up a general point $q$.
The fibers away from $a_0$ are the blow-up of $\P^2$ at $9-n$ points.
The central fiber has two components: W, the blow-up of $\P^2$ at
$10-n$ points and $Z$, the exceptional divisor of the blow-up of $q$.
Over the punctured disk $A^* = A - a_0$ the dual of the relative
dualizing sheaf is a line bundle. One of its flat limits restricts to
the anti-canonical bundle on $W$ and to $\O_{\P^2}(1)$ on $Z$. This
provides us with the desired family.

The limit of the lines are the lines that do not intersect $W \cap Z$.
The conics on the general fiber correspond to lines going through
$p_i$, conics passing through 4 of the points $p_i$ or cubics double
at one $p_i$ and passing through 5 of the other $p_j$. We describe the
limits of conics corresponding to the lines passing through $p_1$. The
others are analogous. In the limit this curve of conics has two
components. One component corresponds to lines passing through
$p_1(0)$ on $W$. The other component consists of the union of the line
$l$ joining $p_1(0)$ and $q$ and a line in $\P^2$ meeting $l$.

The limits of the hyperplane sections have three components. One
component consists of elliptic curves of degree $n$ on $W$ in the
class $3H - \sum_{i=1}^{9-n} E_{p_i}$. One component consists of
conics on $Z$ and rational curves of degree $n-1$ meeting the conics
twice on $W$. The last component corresponds to sections by
hyperplanes that do not contain $W$ or $Z$.

Similar constructions give examples of degenerations of $D_n$ to
$D_{n-k}$ (with various singularities) union a rational cone of degree
$k$. For example, (as B. Hassett pointed out) to obtain a degeneration
of $D_n$ to $D_{n-2}$ with an $A_1$ singularity union a quadric cone
with vertex at the singular point and having a common line with
$D_{n-2}$, blow up the central fiber of $A \times \P^2$ at a general
point, then blow up a general line in the exceptional divisor.  Call
the exceptional divisors of the blow-ups $E_1$ and $E_2$,
respectively. Pick $9-n$ general sections that specialize to the
proper transform of the central fiber and blow them up. Denote the
resulting three-fold by $X$. The linear series $\left| -K_X -2E_1 -E_2
\right|$, where $K_X$ is the canonical bundle of $X$, gives the
desired degeneration.

\medskip

\noindent {\bf VI. Degenerations of $D_n$ to a Veronese surface union
  a rational normal scroll meeting along a conic.} The scroll must
have degree $n-4 < 4$. These scrolls each have at least a one
parameter family of conics (possibly reducible).

Consider a Veronese surface union a scroll $S_{0,1}$, $S_{1,1}$ or
$S_{1,2}$. meeting along a conic. Choose 5, 3, 2 points on their common
conics, respectively. Let $C_1$ be the curve of conics that contain
all but one of the points on the scroll. Let $C_2$ be the curve of
conics that contain the remaining point on the Veronese. This gives a
reducible curve $C= C_1 \cup C_2$ in $X^n$ in the class $d_n$. The
normal bundle $N_{C_i/X^n}$ is generated by global sections. This is
clear for $C_2$ since it lies in a homogeneous locus and follows for
$C_1$ by Lemma \ref{smoothing} after a simple specialization. The
curve can be smoothed to an irreducible curve $\tilde{C}$ in the same
class. By Theorem \ref{classification} the surface $\tilde{S}$
corresponding to $\tilde{C}$ must be $D_n$. Note that the total space
of the family is singular at points that define $C_i$.

\smallskip

\noindent {\bf Alternative description when $n=5$.}  Choose $4$
general sections $p_i(a)$ in $ A \times \P^2$ that specialize to lie
on a line $l$ in the central fiber.  Blow up $A \times \P^2$ along
$l$, then blow up the proper transforms of $p_i(a)$. The general fiber
is the blow up of $\P^2$ at 4 points. The central fiber is $\P^2$
union the blow-up of $F_1$ at 4 points where the two surfaces are
joined along $l$. The dual of the relative dualizing sheaf is a line
bundle away from the central fiber. To obtain a map that extends to
the central fiber we have to twist by the plane in the central fiber.
The limit restricts to $\O_{\P^2}(2)$ on $\P^2$ and to $\O_{F_1}\left(
  e+f-\sum_{i=1}^4 E_i \right)$ on the blow-up of $F_1$. The image is
a Veronese surface union a plane. The total space of the image in $A
\times \P^5$ has 5 singular points. They correspond to the
intersection points of $l$ with the fibers of $F_1$ that are blown
down and with the curve in the class $e+2f$ passing through the $p_i$
which is also blown down. The limits of the lines and conics are
clear. \medskip

\noindent {\bf VII. Degenerations of $D_n$ to the union of two scrolls
  sharing a conic.}  $D_3$'s can specialize to the union of a plane
and a quadric surface. $D_4$'s can specialize to the union of a plane
and $S_{1,2}$ or to the union of 2 quadric surfaces. $D_5$ can
specialize to the union of a plane and $S_{2,2}$ or the union of the
quadric surface and $S_{1,2}$. $D_6$ can specialize to the union of
two $S_{1,2}$, a quadric surface and an $S_{2,2}$ or a plane and
$S_{2,5}$. $D_7$ can specialize to $S_{1,2}$ union $S_{2,2}$, a
quadric surface union $S_{2,3}$, a plane and $S_{2,4}$.  Since further
unbalanced scrolls are limits of balanced scrolls those also arise as
limits.

Take the curve in $X^n$ whose points correspond to incident line pairs
one in each surface. In case one of the scrolls is $\P^2$ take the
lines on $\P^2$ containing a fixed point on the common conic.  This
curve smooths.  When each of the scrolls have a $k_i \geq 1$ parameter
family of conics, the conics on each scroll passing through $k_i -1$
points on the common conic give a different curve which can be
smoothed.

In a general family arising in one of these ways, the limit of
hyperplane sections have three components. Two of the components
correspond to elliptic curves on $S_{k_i,l_i}$ of degree $k_i + l_i
+2$ union $k_j + l_j -2$ fibers on $S_{k_j,l_j}$ meeting the elliptic
curve. The third component corresponds to hyperplane sections by
hyperplanes not containing either of the two components. Since for any
choice of points the curves can be smoothed in $X^n$, we can also
conclude that every elliptic curve of degree $k_i +l_i +2$ occurs as
the limit of some family.  \medskip

\noindent {\bf $D_8$.} The arguments for the degenerations in I, III
and V apply to $D_8$ verbatim. However, since the anti-canonical
embedding of $\P^1 \times \P^1$ is also a smooth degree 8 surface in
$\P^8$ that contains a curve of conics in the class $d_8$, the
arguments in II, VI and VII only show that one of the two surfaces
degenerates to the potential limit. Ciro Ciliberto pointed out to
us that it is possible to modify the alternative construction in II
using a Cremona transformation $(3,3,3,1)$ to obtain a degeneration of
$D_8 $ to a scroll with a double line. In VI when we smooth the
Veronese union $S_{2,2}$ we obtain $D_8$ since the limit surface does
not contain two rulings by conics where a conic from one ruling meets
every conic in the other ruling as a limit of $\P^1 \times \P^1$
should. In VII we note that the union of two $S_{2,2}$ does not smooth
to $D_8$ since it does not contain a line meeting every conic as a
limit of $D_8$ should.  \smallskip

\noindent {\bf Degenerations of the Veronese surface.} Since the
Veronese surface appears as a component of the limits of $D_n$, we
mention the non-degenerate two component limits of it. They are the
union of a plane and a cubic scroll where the cubic scroll meets the
plane along the directrix or the union of two quadric cones that share
a vertex and a common line. Both cases occur. The Veronese surface
degenerates to a cone over a rational normal quartic. The latter is a
further specialization. To obtain the former limit carry out the usual
construction by blowing up $A \times \P^2$ at a point on the central
fiber.

Since the surface spans $\P^5$ the components must be scrolls meeting
along a line. They can have degrees 2, 2 or 1, 3. Since a cone over a
twisted cubic is a limit of $S_{1,2}$, we can assume that the cubic
surface is smooth. The cubic plane pair cannot be joined along a fiber
line and the quadrics have to be both singular with a common vertex.
The former cannot happen because any two such surfaces are
projectively equivalent. The dimension of the locus of pairs of a
cubic surface union a plane meeting it along a fiber is too large. To
show that the quadrics have to be as described we only need to show
that the union of two quadric cones that share a common line but have
distinct vertices cannot be a limit of Veronese surfaces.

Each Veronese surface has a two-parameter family of conics, hence gives
rise to a $\P^2$ in $X^5$. The flat limit of the conics has to be a
surface in $X^5$ connected in codimension 1. If a pair of quadric
cones with distinct vertices were a limit, then the two surfaces would
need to share a curve of conics (possibly reducible). Since this is
not the case, we conclude that the vertices have to coincide.

A forthcoming paper of Dan Avritzer should elucidate the enumerative
geometry of the Veronese surface. We are thankful to him for
conversations on the degenerations of the Veronese surface.

\section{Examples of Counting Del Pezzo Surfaces}

In this section we illustrate with a few examples how to use our
knowledge of the degenerations of Del Pezzo surfaces to study their
enumerative geometry.  \medskip

\noindent {\bf Example 1: Counting cubic surfaces in $\P^4$.} By Lemma
\ref{dimdelpezzo} the dimension of the space of cubic surfaces in
$\P^4$ is 23. We can ask for the number of $D_3$'s containing $r$
points and meeting $23-2r$ lines. \smallskip

\begin{figure}[htbp]
\begin{center}
\epsfig{figure=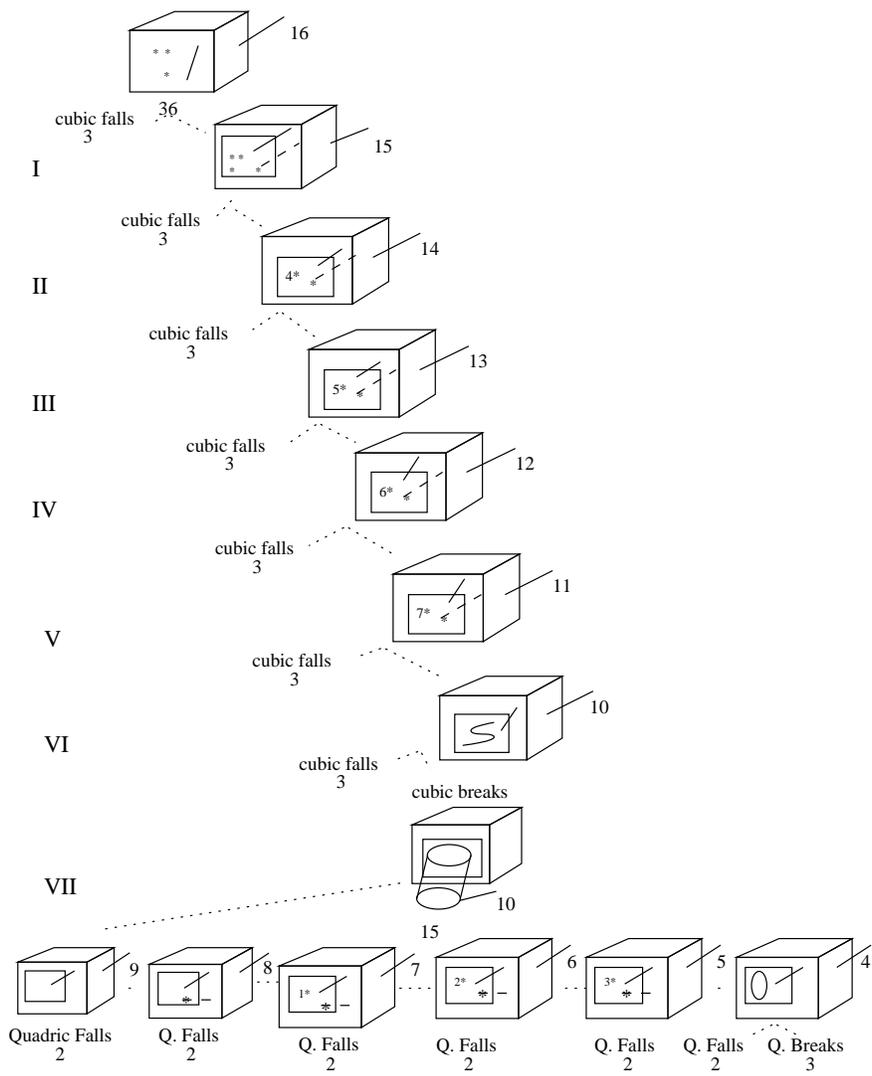}
\end{center}
\caption{Cubic surfaces in $\P^4$ containing 3 points and meeting 17
  lines (Example 1) }
\label{Figure 3}
\end{figure}

\noindent {\bf The number of cubic surfaces containing 3 points and meeting 17
  lines.} We outline how to see that there are 36 cubic surfaces
satisfying the required incidences using degenerations. (See Figure
\ref{Figure 3}.) Fix a hyperplane $H$ in $\P^4$. \smallskip

\noindent {\bf Step I.} Specialize the three points and a line $l_1$ to $H$.

$\bullet$ Some cubic surfaces can lie in $H$.  These surfaces must meet
the 16 intersection points of $H$ with the lines outside $H$. They
must also contain the original 3 points. There is a unique $D_3$
containing the 19 points. It counts with multiplicity 3 for the choice
of intersection point of the surface with $l_1$.

$\bullet$ If a cubic surface does not lie in $H$, then its hyperplane
section must be contained in the plane $P$ spanned by the three points
in $H$, so the surface must contain $P \cap l_1$. 

\noindent {\bf Step II.} Specialize a second line $l_2$ to $H$.
Some cubics can now lie in $H$. Such a cubic must meet 19 points---the 4
points in the plane $P$ and the 15 points of intersection of $H$ with
the lines outside $H$. These points impose independent conditions on
cubics, so there is a unique solution counted with multiplicity 3 for
the choice of intersection with $l_2$.  \smallskip

\noindent {\bf Step VI.} This pattern continues until we specialize 6
lines to $H$. After we specialize 6 lines, if the cubic does not lie
in $H$, then the hyperplane section has to be the unique cubic curve
in $P$ passing through the 9 points in $P$. When we specialize the
next line, either the cubic lies in $H$ or it must break into $P$ and
a quadric meeting the rest of the lines. By a dimension count this is
the first stage where reducible solutions occur.  \smallskip

\noindent {\bf Step VII.} We are reduced to counting quadric surfaces
meeting 10 lines and containing a conic in common with $P$. Further
degeneration shows that there are 15 such quadrics. Briefly,
specialize a line $l_8$ to $H$.  Either the quadric lies in $H$ or it
contains the point of intersection of $l_8$ with $P$. If the quadric
lies in $H$, then it must contain the intersection points of the 9
remaining lines with $H$. This uniquely determines the quadric. It
counts with multiplicity two for the choice of intersection with
$l_8$.

The same pattern continues until we have 5 lines remaining outside.
Then the hyperplane section of the quadric is the unique conic passing
through the 5 points in $P$. Once we specialize another line, the
quadric has to either lie in $H$ or break into a union of $P$ and
another plane having a common line with $P$ and meeting the remaining
4 lines. The latter number is 3 by elementary Schubert calculus.  We
conclude that there are 36 cubic surfaces in $\P^4$ meeting 17 lines
and containing 3 points. We will later verify the multiplicities (see \S 5.2).

\medskip

\noindent {\bf Example 2: Counting $D_4$'s in $\P^4$.} By Lemma
\ref{dimdelpezzo} the dimension of the locus of $D_4$'s in $\P^4$ is
26. We can ask for the number of $D_4$'s containing $r$ general points
and meeting $26-2r$ general lines.  \smallskip

\begin{figure}[htbp]
\begin{center}
\epsfig{figure=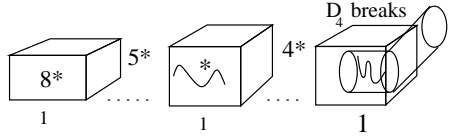}
\end{center}
\caption{One $D_4$ containing 13 general points (Example 2A)}
\label{Figure 4}
\end{figure}

\noindent {\bf A. The number of $D_4$'s in $\P^4$ containing 13 points.} 
Specialize the points to a hyperplane $H$ of $\P^4$. No reducible
surfaces satisfy all the incidences until we specialize 9 points to
$H$. After we specialize 8 points to $H$, the hyperplane section of a
solution must be the unique elliptic quartic containing the 8 points.
When we specialize a ninth point to $H$, Bezout's Theorem forces the
surface to break into 2 quadric surfaces. The quadric surface $Q$ in
$H$ is determined. There is a unique quadric in the $\P^3$ spanned by
the 4 points not in $H$ containing the intersection of $Q$ with the
$\P^3$ and the 4 points. After we verify the multiplicity claims, we
can conclude that there is a unique $D_4$ containing 13 general
points.  \smallskip

\begin{figure}[htbp]
\begin{center}
\epsfig{figure=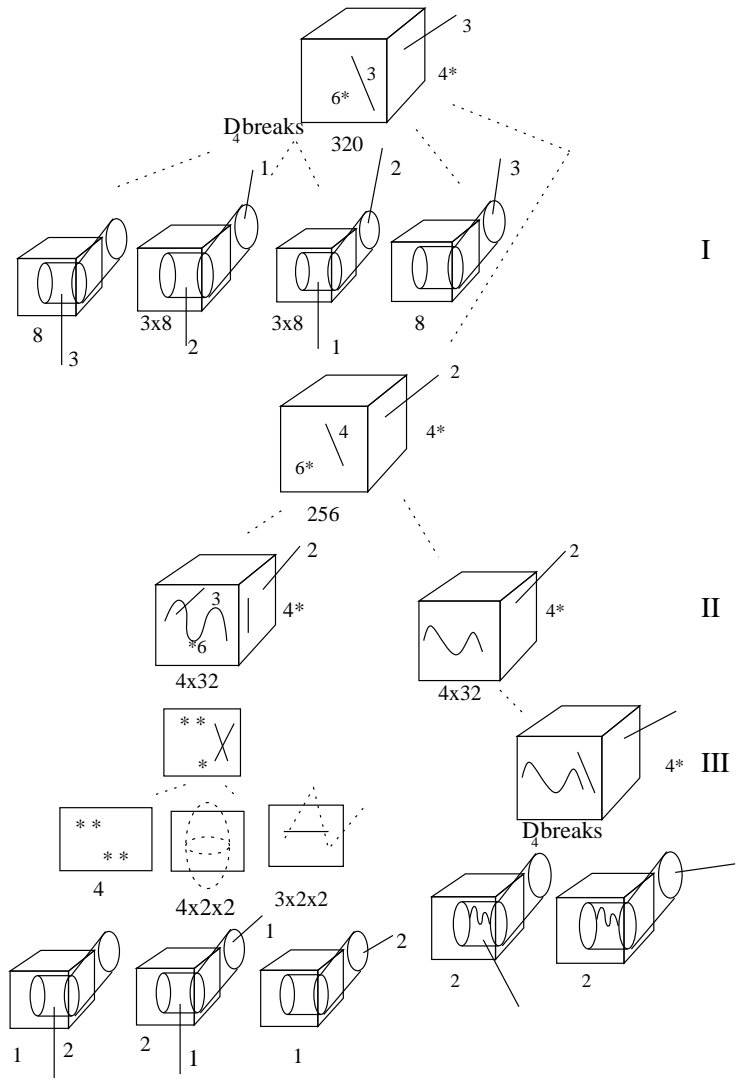}
\end{center}
\caption{Counting $D_4$ surfaces containing 10 general points and
  meeting 6 general lines (Example 2B)}
\label{Figure 5}
\end{figure}

\noindent {\bf B. The number of $D_4$'s in $\P^4$ containing 10 points and
  meeting 6 lines.} (See Figure \ref{Figure 5}.) Fix a hyperplane $H$
of $\P^4$.

\smallskip 

\noindent {\bf Step I.} Specialize 6 points and 3 lines $l_1, \  l_2,
\ l_3$ to $H$. This is the first stage where reducible solutions
occur: there can be the union of two quadric surfaces meeting along a
conic.  Of the three lines outside $H$, $3,\ 2, \ 1$ or 0 of them
might meet the quadric in $H$. The remaining $0, \ 1, \ 2$ or 3 lines
need to meet the quadric outside $H$. In each case there is a
multiplicity of $8$ for the choice of intersection points of the lines
in $H$ with the quadric in $H$. In two cases there is a combinatorial
choice of 3 for which of the lines meet the quadric in $H$. The
surfaces are uniquely determined. \smallskip

\noindent {\bf Step II.} If a solution is still irreducible, then its
hyperplane section in $H$ is an elliptic quartic curve $C$ meeting the
6 points and the lines $l_1, l_2, l_3$. Specialize a line $l_4$ to
$H$. The new reducible solutions must contain a curve $C$ as
described. We specialize three points and two lines $l_1, l_2$ to lie
in a plane $P$ in $H$. Either $C$ passes through the intersection
point of $l_1$ and $l_2$ or it must have a component in $H$. This
component can either be a conic or a line. The residual component must
be a conic or a twisted cubic, respectively. The number of these can
be determined using the algorithm in \S 7 of \cite{vakil:rationalelliptic}.
Finally, of the lines outside $H$, $2,\ 1$ or 0 of them might meet the
quadric in $H$.  Considering all the cases we see that there are 128
reducible surfaces at this stage. \smallskip

\noindent {\bf Step III.} If a solution is still irreducible, then
its hyperplane section in $H$ must be one of the 32 elliptic quartic
curves containing 6 points and meeting 4 lines
(\cite{vakil:rationalelliptic} \S 8.3). We specialize a fifth line
$l_5$ to lie in $H$. The Del Pezzo surfaces now have to break into a
union of two quadrics. The sixth line can either meet the quadric in
$H$ or the quadric not lying in $H$. In each case there is a unique
surface, appearing with multiplicity 2 for the choice of intersection
of $l_5$ with the quadric in $H$. We conclude that there are 320
quartic Del Pezzo surfaces in $\P^4$ containing 10 general points and
meeting 6 lines. \medskip

\begin{figure}[htbp]
\begin{center}
\epsfig{figure=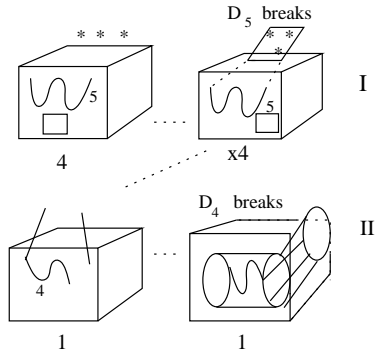}
\end{center}
\caption{Counting $D_5$'s containing a quintic elliptic curve, 3 points and
  meeting a plane (Example 3)}
\label{Figure 6}
\end{figure}

\noindent {\bf Example 3: Counting $D_5$'s in $\P^5$ containing an
  elliptic quintic curve.} As a final illustration of the type of
enumerative problems one can hope to answer using degenerations, we
find the number of $D_5$'s in $\P^5$ containing an elliptic quintic
curve, three points and meeting a plane. (See Figure \ref{Figure 6}.)

\smallskip

\noindent {\bf Step I.} Specialize the plane $P$ to the hyperplane $H$
spanned by the elliptic quintic. The surface breaks into a union of
$D_4$ and a plane $\Pi$. (By a dimension count it cannot break into a
cubic scroll union a quadric surface. The other possibilities in
Proposition \ref{junk} are either further degenerations of $D_4$,
hence are excluded by a dimension count or do not contain any elliptic
quintic curves---e.g. a Veronese surface or a quartic scroll.) We
reduce the problem to counting $D_4$'s containing an elliptic quintic
$C$ and a disjoint line $l$ ($l= H \cap \Pi$). We get a multiplicity
of 4 for the choice of intersection of $P$ with any limit $D_4$.
\smallskip

\noindent {\bf Step II.} Specialize the elliptic quintic to the union
of an elliptic quartic and a general line. $D_4$ must become reducible
by Bezout's Theorem since the hyperplane spanned by the elliptic
quartic meets $l$. The surface must break into a union of quadrics and
they are both uniquely determined. We conclude that there are four
quintic Del Pezzo surfaces containing an elliptic quintic, three
points and meeting a plane.  \smallskip

\noindent {\bf A non-example.} Even when we impose only point
conditions on surfaces $D_n$ for $n > 4$, at each stage new and more
complicated degenerations appear. One can make arbitrarily long lists
of degenerations to count more cases until the dimension counts or the
multiplicity calculations lose their rigor. We will instead limit
ourselves to simple enumerative problems. We give, however, an example
to illustrate the complications due to curve conditions.

\begin{figure}[htbp]
\begin{center}
\epsfig{figure=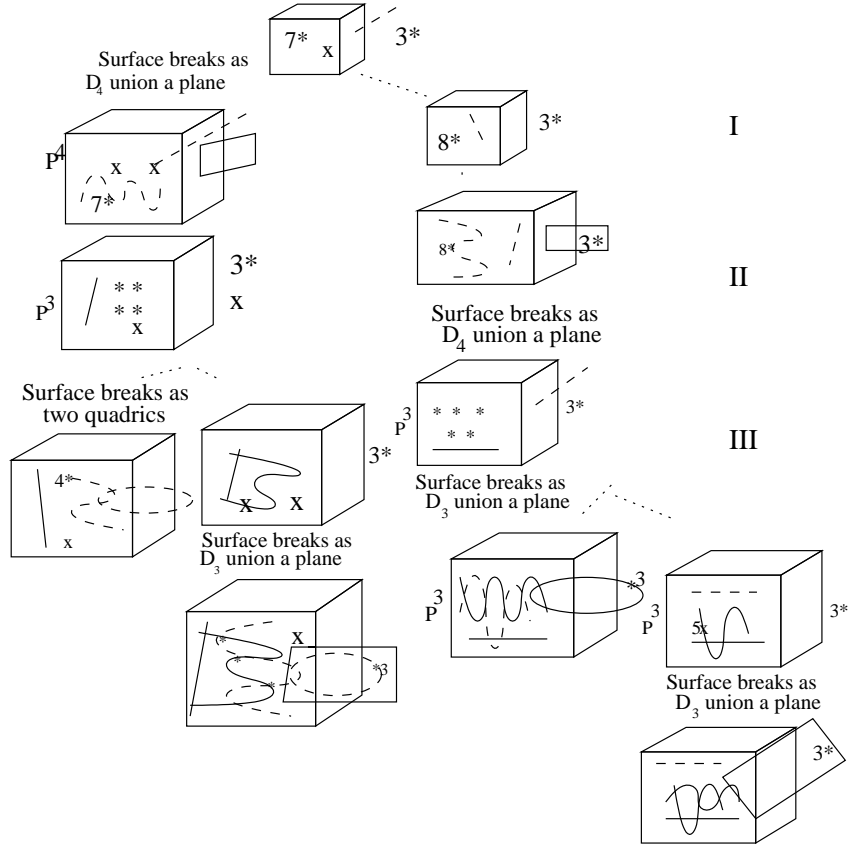}
\end{center}
\caption{Degenerations of  $D_5$'s containing 11 points and
  meeting a line}
\label{Figure 7}
\end{figure}

\noindent {\bf The degenerations of $D_5$ in $\P^5$ that contain
  11 points and meet a line $l$.} (See Figure \ref{Figure 7}.) When we
specialize 8 points to $H$, the surface can break into a plane $P$
union a $D_4$. The $D_4$ must contain an elliptic quintic curve (the
limit of the hyperplane sections of $D_5$) passing through the first 7
points (marked by *), the line $l_P = P \cap H$ and 2 other points
(marked by x)---the 8th point and the point $l\cap H$. We specialize
$l_P$, 4 of the 7 points (*) and one of the points (x) to a $\P^3$. If
the surface does not break, then the hyperplane section is $l_P$ union
the twisted cubic $C$ meeting it twice and containing 5 points (4 *
and one x). When we specialize the last point (x), $D_4$ breaks into a
cubic surface union a plane $P'$. The elliptic quintic must break into
a union of a twisted cubic $C'$ with a conic outside $\P^3$ meeting it
twice since $3$ of the points (*) are still outside the $\P^3$. Now we
need to determine cubic surfaces containing the line $l$ with the
twisted cubic $C$ meeting it twice, containing a line $l'$ which is
the intersection of $P'$ with the $\P^3$ and in addition containing a
twisted cubic $C'$ which meets $l'$ twice and contains the 4 points
(*) on $C$. The other cases are similar. As this example indicates
when $n$ gets larger, the curve conditions on $D_n$ get very
complicated making it very hard to continue the degenerations.

\section{The Enumerative Geometry of $D_n$}

In this section we carry out the dimension and multiplicity
calculations necessary to justify calculations similar to ones in \S
4. 

\subsection{Dimension Counts}

\noindent {\bf The building blocks.} We calculate the dimension of
relevant loci in the Hilbert scheme of surfaces.

\begin{lemma}\label{conicH}
  The dimension $D$ of the locus of surfaces $S$ in $\P^N$ containing an
  irreducible, reduced conic in a fixed hyperplane $H$ is as follows:

1. If $S$ is a scroll $S_{2,l}$, $l \geq 2$, then $D=N(l+4) -3$.

2. If $S$ is a scroll $S_{1,l}$, $1 \leq l \leq 2$, then $D= N(l+4) -
3$.

3. If $S$ is the Veronese surface, then $D=6N-4$.
\end{lemma}

\begin{lemma} \label{lineH}
  The dimension $D$ of the locus of surfaces $S$ in $\P^N$ containing
  a line $l$ in a fixed hyperplane $H$ is as follows:
  
1. If $S$ is a rational cone $S_{0,l}$, then $D = N(l+2) - 5$.
  
2. If $S$ is a smooth Del Pezzo surface $D_n$, then $D= N(n+1) - n +
8$.
   
3. If $S$ is a cone over an elliptic normal curve of degree $k$
tangent to $H$ everywhere along $l$, then $D=N(k+1) -2$.
\end{lemma}

\begin{lemma}\label{inside}
  The dimension $D$ of the locus of pairs $(S,C)$ in $\P^N$ where $S$
  is a surface and $C$ is a curve on it is as follows:

1. When $S$ is a Del Pezzo surface $D_n$ and $C$ is an elliptic curve
of degree $n$ (resp. $n+1$) , then $D= N(n+1) + 10$ (resp. $D=N(n+1) +
11$),
   
2. When $S$ is a scroll $S_{k,l}$, $l-k \leq 2$, and $C$ is an
elliptic curve in a bisection class $2e + (l-k+2) f$, then $D=
N(k+l+2) + 2k+4 - \delta_{k,l}$
   
3. When $S$ is a Del Pezzo surface $D_n$ and $C$ is a rational curve
of degree $n-1$, then $D= N(n+1) + 9$.
\end{lemma}

\noindent {\bf Proof:} To prove the lemmas consider maps from $\P^2$,
$F_{l-2}$, $F_{l-1}$ or the blow-up of $\P^2$ at $9-n$ points to
$\P^N$, up to isomorphism, given by the linear series
$\left| \O_{\P^2}(2) \right|$, $\left| \O_{F_{l-2}}\left( e+lf \right)
  \right|$, $\left| \O_{F_{l-1}} \left( e+lf \right) \right|$,
$|-K|$ respectively. In the cases where the surface is required to
contain a curve in a fixed hyperplane $H$, we assume that the map is
given by $(s_0, \cdots , s_N)$, where $s_0$ corresponds to $H$. The
lemmas follow from the cohomology calculations in \S 2 and \S2 of
\cite{coskun1:degenerations}. $\Box$

 \smallskip

\noindent {\bf Gluing.} We now prove that gluing surfaces along
projectively equivalent rational curves to form a
tree imposes the expected number of conditions. \smallskip

\noindent {\bf Notation.} For a variety $X \in \P^N$, let
$\mathcal{H}_X$ denote the locus in the Hilbert scheme parametrizing
varieties projectively equivalent to $X$. Let $\mbox{Rat}_d (X)$ be an
irreducible subscheme of the scheme of rational normal curves of
degree $d$ on $X$. For a variety $Z$ such that $[Z] \in \mathcal{H}_X$
let $\mbox{Rat}_d (Z)$ denote the transform of $\mbox{Rat}_d (X)$
under the projective linear transformation that takes $Z$ to $X$. Let
$(p_i^X)_{i=1}^s$, $0 \leq s \leq 3$, denote distinct points on a
rational normal curve.

\begin{lemma}\label{glue}
  Let $X_1$ and $X_2$ be two varieties in $\P^N$. In the incidence
  correspondence $I := $
\begin{equation*}
\begin{split}
\left\{ \left( Z_1,C_{Z_1}, (p_i^{Z_1})_{i=1}^s, Z_2, C_{Z_2},
(p_i^{Z_2})_{i=1}^s  \right) : Z_j \in
\mathcal{H}_{X_j}, C_{Z_j} \in \mathrm{Rat}_d(X_j),  p_i^{Z_j} \in
C_{Z_j} \right\}
\end{split} 
\end{equation*}
the locus $C_{Z_1} = C_{Z_2}$ and $p_i^{Z_1} = p_i^{Z_2}$ for all $1
\leq i \leq s$ has codimension $$N(d+1) + d-3 +s.$$    
\end{lemma}  

\noindent {\bf Proof.} Without loss of generality we can assume that
$s=0$. $I$ maps to $\mbox{Rat}_d(\P^N) \times \mbox{Rat}_d(\P^N)$ by
projection. The fibers are equivalent under the diagonal action of
$\P GL(N+1)$. Since the locus of interest is the inverse image of the
diagonal, the lemma follows. $\Box$ \smallskip

\noindent {\bf Notation.} Let $H$ and $\Pi$ denote two hyperplanes in
$\P^N$. Let $\Sigma_{a_j}^j$ and $\Omega_{b_i}^i$ be collections of
general linear subspaces of $H$ and $\P^N$ of dimension $a_j$ and $b_i$,
respectively. Similarly, let $\Lambda^{j}_{a_j}$ and
$\Gamma_{a_{j'}}^{j'}$ be collections of general  linear spaces of $\P^N$ and
$\Pi$, respectively.  We will usually omit the dimension from the
notation.  We denote connected curves of arithmetic genus 1 by $E$ and
connected curves of arithemetic genus 0 by $R$. To denote their degree
we append a number in parentheses.

Let $\mathcal{H}(\P^N, D_n)$ be the component of the Hilbert scheme
whose general point corresponds to a smooth Del Pezzo surface $D_n$.
Let $\mathcal{E}(\P^N, m)$ denote the component of the Hilbert scheme
whose general point represents a smooth elliptic curve of degree $m$
in $\P^N$. Let $\mathcal{HE}(\P^N,D_n, m)$ be the incidence
correspondence of pairs
$$\left\{ \left([D_n], [E(m)] \right) \in \mathcal{H}(\P^N,D_n) \times
\mathcal{E}(\P^N,m) : E \subset D_n \right\}$$
where the elliptic curve $E$ is a closed subscheme of the Del Pezzo
surface $D_n$. Finally, let $\mathcal{U}(\P^N,D_n,m, I,J)$ denote the
$(I,J)$ pointed universal surface curve pair over
$\mathcal{HE}(\P^N,D_n, m)$ defined by
$$\left\{ \left(D_n, E , (q_i)_{i=1}^I, (p_j)_{j=1}^J\right): p_j \in
  E \subset D_n, q_i \in D_n, (D_n,E) \in \mathcal{HE}(\P^N,D_n, m)
\right\} $$
where $p_j$ and $q_i$ are points of the curve $E$ and
surface $D_n$, respectively.

\smallskip

\noindent {\bf The space of Del Pezzo surfaces.} Let $\mathcal{D}_n
(\P^N, I_1,I_2,J)$ denote the closure in $\mathcal{U} (\P^N,D_n,n,I_1
+ I_2,J)$ of
$$
\left\{ \left( D_n,E, (q_i)_{i=1}^I, (p_j)_{j=1}^J \right): E = D_n
  \cap H, q_i \subset \Omega^i, p_j \subset \Sigma^j \subset
  H\right\}$$
where $D_n$ is a smooth Del Pezzo surface, $E$ is its
hyperplane section in $H$ and the marked points are requied to lie in
the designated linear spaces. This notation is bad since many
different possibilities are denoted by the same symbol. Since in any
given enumerative problem the dimensions of the linear spaces will be
predetermined we will use it as a shorthand.  \smallskip

Let $I_1 \cup I_2 =I$ and $J_1 \cup J_2 =J$ be two partitions.

\noindent {\bf The space of scroll pairs.} Let $\mathcal{S} (\P^N,
k_1,l_1,k_2,l_2,I_1,I_2, J_1,J_2)$ denote the closure in $\mathcal{U}
(\P^N,D_n,n,I_1+I_2,J_1 + J_2)$ of the locus
\begin{equation*}
\begin{split}
 \{ & \left( S_{k_1,l_1} \cup S_{k_2,l_2},
      E \left( k_1+l_1+2 \right) \cup F_1 \cup \cdots \cup F_{k_2+l_2-2},
      q_i, p_j \right) : S_{k_1,l_1} \subset H \\ &
    E \left( k_1+l_1+2 \right) \subset S_{k_1,l_1}, F_1 \cup \cdots \cup
    F_{k_2+l_2-2} \subset S_{k_2,l_2} \cap H, S_{k_1,l_1} \cap
    S_{k_2,l_2} = R(2) \\ & q_i \in \Omega^i \cap S_{k_r,l_r} \ 
    \mbox{for} \ i \in I_r, \ p_j \in\Sigma^j \cap
    E \left( k_1+l_1+2 \right) \ \mbox{for} \ j \in J_1, \\ & p_j
    \in\Sigma^j \cap \left( F_1 \cup \cdots \cup F_{k_2+l_2-2} \right) \ 
    \mbox{for} \ j \in J_2 \}
\end{split}
\end{equation*}
of pairs of scrolls meeting along a conic $R(2)$, where the scroll
$S_{k_1,l_1}$ is in $H$ and contains an elliptic curve $E(k_1+l_1+2)$
and the other scroll is outside $H$ and its intersection with $H$
consists of the conic $R(2)$ and the fibers $F_1, \cdots,
F_{k_2+l_2-2}$ and the marked points lie in the designated linear
spaces. The elliptic curve $E$, needless to say, meets all the fibers
$F_i$. We also assume that $k_2+l_2 > 1$ to ensure that $S_{k_2,l_2}$
can lie outside $H$.  \smallskip

\noindent {\bf The space of $D_{n-1}$ union a plane.} Let $\mathcal{PD}_{n-1}
(\P^N, I_1,I_2,J)$ denote the closure in $\mathcal{U}          
(\P^N,D_n,n,I_1+I_2,J)$ of the locus
\begin{equation*}
\begin{split}
  \{ &(D_{n-1} \cup \P^2, E(n), q_i, p_j) :D_{n-1} \cap \P^2 = R(1),
  E(n) \subset D_{n-1} \subset H, \\ & q_i \in \Omega^i \cap D_{n-1} \ 
  \mbox{for} \ i \in I_1, q_i \in \Omega^i \cap \P^2 \ \mbox{for} \ i
  \in I_2, p_j \in\Sigma^j \cap E(n)\}
\end{split} 
\end{equation*}
where $D_{n-1}$ is a smooth Del Pezzo surface in $H$, $E(n)$ is an
elliptic curve on $D_{n-1}$ and the marked points lie in the
designated linear spaces.

Let $\mathcal{D}_{n-1} \mathcal{P} (\P^N, I_1,I_2,J_1,J_2)$ be the
variant where $\P^2$ is in $H$ and $D_{n-1}$ is outside, the two meet
in a line $l$ and the marked curve is a conic in $\P^2$ meeting the
hyperplane section of $D_{n-1}$ residual to $l$ in two points.
\smallskip

\noindent {\bf To generalize} the discussion below to more cases one
has to formulate similar loci corresponding to other reducible
surfaces (possibly with more components) that occur in \S 3 like pairs
of a Veronese and a scroll or an elliptic cone and a rational cone. In
addition one has to allow for the surfaces outside to have tangencies
with $H$ along their common curves with the surfaces in $H$ and the
``limit curve'' in the components in $H$ to have correspondingly
larger degree. \smallskip

\noindent {\bf The space of marked $D_n$.} Let $D_n (\P^N,m, I, J, J')$
denote the closure in $\mathcal{U}(\P^N,D_n,m, I, J+J')$ of the locus

$$\left\{ \left( D_n, E, (q_i)_{i=1}^I, (p_{j})_{j=1}^{J},
    (o_{j'})_{j'=1}^{J'} \right) : q_i \in D_n \cap \Omega^i , p_j \in
  \Lambda^j \cap E, o_{j'} \in \Gamma^{j'} \cap E \right\} $$
where
$D_n$ is a smooth Del Pezzo surface, $E$ is a degree $m$ elliptic
curve on it and the marked points lie in the designated linear spaces.
The indices $i$ are reserved for points on the surface, the indices
$j$ indicate points on $E$, but not in $\Pi$. Finally, the indices
$j'$ designate points that lie both on $E$ and in $\Pi$.

Let $D_n (\P^N, m, I, J,J', \O(1))$ denote the analogous space, but
where in addition the points $o_{j'}$ satisfy $\sum_{j'=1}^{J'} o_{j'}
= \O_{E}(1)$ in the Picard group of the elliptic curve $E$. 

Let $S_{k,l} ( \P^N, k+l+2, I, J, J')$ and $S_{k,l} ( \P^N,
k+l+2, I, J, J', \O(1))$ denote the analogous space where $D_n$ is
replaced by a scroll $S_{k,l}$ and the elliptic curve has degree
$k+l+2$.

Let $D_n (\P^N, r_1 + r_2 = m , I,J_1,J_2,J_1',J_2')$ denote the
closure of the locus where $E$ is a pair of rational curves of degrees
$r_1$ and $r_2$ meeting at two points in $D_n (\P^N, m, I,
J_1+J_2,J_1'+J_2')$ and the conditions are distributed between the
rational curves according to a partition.

We define the analogous locus $S_{k,l}(\P^N; r_1+r_2=k+l+2;
I,J_1,J_2,J_1', J_2')$ for scrolls.

 \smallskip

\noindent {\bf The divisors.} The space $\mathcal{D}_n
(\P^N, I,J)$ has a natural Cartier divisor
$$
D_H (\P^N,I,J) := \{ (D_n, E, q_i, p_j) \in \mathcal{D}_n (\P^N, I,J): q_I \in
H \}$$
defined by requiring one of the marked points on the surface to
lie in $H$.

Let $D_{\Pi}(\P^N,D_n,I,J,J')$ on $D_n (\P^N, n+1, I, J, J')$ and
$D_{\Pi} (\P^N, S_{k,l},I,J,J')$ on $S_{k,l} ( \P^N, E(k+l+2), I, J,
J')$ defined by letting $p_J$ lie in $\Pi$ be analogous Cartier
divisors.

There is a natural map $\phi$ from $\mathcal{U}(\P^N,D_n,m,I,J)$ to
$\mathcal{H}(\P^N,D-n)$ given by projection.

\begin{definition}
  A divisor $D$ of a subscheme $A$ in $\mathcal{U}(\P^N,D_n,m,I,J)$ is
  called {\bf enumeratively relevant} if the image of $D$ under $\phi$
  has codimension 1 in the image of $A$.
\end{definition}
One can list the components of $D_H(\P^N,I,J)$ that are enumeratively
relevant Weil divisors whose general point corresponds to a surface
that has reduced components and spans $\P^n$. Rather than give long
lists, we will demonstrate how one produces such lists in examples
similar to the ones in \S 4 and indicate what other type of behavior
can occur when we require fewer point conditions.

The following lemma, which is a consequence of Kleiman's
Transversality Theorem \cite{kleiman:transverse} (see Proposition 6.1
in \cite{coskun1:degenerations}) allows us to carry out the dimension
calculations when $I=1$, $\Omega^j = H$ and $ \Sigma^1 = \P^N$.

\begin{lemma}\label{reduction}
Let $A$ be a reduced, irreducible subscheme of
$\mathcal{D}_n(\P^N,I,J)$ and let $p$ be one of the labeled
points. Then there exists a Zariski open subset $U$ of the dual
projective space $\P^{N*}$ such that for all hyperplanes $[H] \in U$,
the intersection $A \cap \{ p \in H \}$ is either empty or reduced of
dimension $\dim A -1$. 
\end{lemma}

\noindent We describe the components of $D_H(\P^N,I,J)$ in
$\mathcal{D}_n (\P^N,I,J)$ relevant to our examples.

$$1. \ \mathcal{D}_n (\P^N, I- 1, J+1) $$
where $ \Sigma^{J+1} :=
\Omega^I$ and the rest of the data is identical. To check that this
locus is a divisor we use the above reduction. In that case the
dimension for surfaces does not change; however, the moduli for the
points is one less. Since surfaces which are limits of $D_n$ which do
not contain a component in $H$ lie in the closure of $1$, we can now
assume that all other limit surfaces have a component in $H$.

$$2.\ \mathcal{D}_n (\P^{N-1}, I, J \geq n)$$
where $ q_I \in
\Omega^I$ and $q_i \in H \cap \Omega^i, \ i < I $. Since the dimension
of $D_n$ in $\P^{N-1}$ is $n^2+10+(n+1)(N-1-n)$, the choice for
hyperplane section has dimension $n$ and the points have the same
moduli, this locus is a divisor. It is enumeratively relevant when $J
\geq n$. When $J<n$, the hyperplane sections move in a positive
dimensional linear series. The image of $\phi$ has positive
dimensional fibers. We can now assume that the other components of
$D_H$ consist of reducible surfaces at least one component of which
lies in $H$ and one component lies outside $H$. \smallskip

$$3. \ \mathcal{S} (\P^N, k_1,l_1,k_2,l_2, I_1,I_2, J_1,J_2), \ J_1+
J_2 \geq 8 - (k_2 + l_2 -2) + \delta^{(0,1)}_{(k_1,l_1)}, \ l_1 - k_1
\leq 2$$
where $\sum_i (k_i + l_i) =n, \ q_I \in S_{k_1,l_1}$ and the
rest of the marked points are distributed between the two components
according to some partition. By Lemmas \ref{conicH}, \ref{inside} and
\ref{glue} and the constructions in \S 3 the locus 3 is a divisor. It
is enumeratively relevant when $J \geq 8 - (k_2 + l_2-2)$ unless
$(k_1,l_1)=(0,1)$ in which case it is enumeratively relevant when $J
\geq 9 - (k_2 + l_2-2)$ . \smallskip

$$4. \ \mathcal{PD}_{n-1}(\P^N, I_1,I_2,J \geq n > 3)$$
where $q_I \in
D_{n-1}$ and the other points are distributed among the two components
according to some partition. By Lemmas \ref{lineH}, \ref{inside} and
\ref{glue} it is a divisor. It is enumeratively relevant when $J \geq
n$.

\subsection{The ``Algorithm'' for Counting $D_3$ and $D_4$}
Now we outline the argument for counting $D_3$ and $D_4$ in greater
detail. We begin with some lemmas necessary for multiplicity
calculations.

\begin{lemma}\label{counting}
Let $E$ be an elliptic curve of degree  $n+1$ on $D_n$ in the
class $-K + E_i$ then 
\begin{enumerate}
\item $H^1(D_n , N_{D_n/\P^N})=0$,

\item $H^1(D_n , N_{D_n/\P^N}(-E))=0$

\end{enumerate}
\end{lemma}

\noindent {\bf Proof:} Suppose $H^i (D_n, T_{\P^N} \otimes \O_{D_n}) =
0$ and $H^i (D_n, T_{\P^N} \otimes \O_{D_n} (-E)) =
0$  for $i > 0$. Then using the standard short exact sequence
$$0 \rightarrow T_{D_n} \rightarrow T_{\P^N} \otimes \O_{D_n}
\rightarrow N_{D_n/P^N} \rightarrow 0$$
we conclude that $H^1(D_n ,
N_{D_n/\P^N}) \cong H^2(D_n, T_{D_n})$. The analogous statements hold
when we twist the sheaves by $-E$.  Tensoring the Euler sequence 
$$
0 \rightarrow \O_{\P^N} \rightarrow (\O_{\P^N}(1))^{N+1}
\rightarrow T_{\P^N} \rightarrow 0$$
by $\O_{D_n}$ and $\O_{D_n}(-E)$
we conclude that $$H^i (D_n, T_{\P^N} \otimes \O_{D_n}) = 0 \ \ \
\mbox{and} \ \ H^i
(D_n, T_{\P^N} \otimes \O_{D_n} (-E)) = 0$$ for $i > 0$. Here we are
using the fact that $H^2(D_n, \O_{D_n} (-E)) = 0$. This follows by
Serre duality.

Finally, by Serre duality $H^2(D_n, T_{D_n}) = 0$ and $H^2(D_n,
T_{D_n}(-E)) = 0$. The lemma follows. $\Box$ \smallskip

\begin{lemma}\label{counting3}
  Let $D_3$ be a smooth cubic surface. Let $E$ be a curve in the class
  $-K$ on $D_3$. Then $H^1(D_3, N_{D_3, \P^N} (-E)) = 0$
\end{lemma}

\noindent {\bf Proof:} $N_{D_3/\P^N} \cong \O_{D_3}(-3K) \oplus
(\O_{D_3} (-K))^{N-3}$. We twist the normal bundle by $K$. Since $D_3$
is a rational surface, $h^1(\O_{D_3}) = 0$; and $h^1(\O_{D_3}(-2K))$
vanishes by the Kodaira Vanishing Theorem. $\Box$ \smallskip

\begin{lemma} \label{counting4} 
  Let $Q$ be a smooth quadric surface. Let $E$ be a curve in the class
  $\O_Q(2,2)$ on $Q$. Then $H^1(Q, N_{Q, \P^N} (-E)) = 0$
\end{lemma}

\noindent{\bf Proof:} Since $N_{Q/\P^N} \cong \O_Q (2,2) \oplus (\O_Q
(1,1))^{N-3}$ the lemma follows from the cohomology of $Q$
(\cite{hartshorne:book} Ex. III.5.6). $\Box$ \bigskip

\noindent {\bf Counting $D_3$.} We now analyze the case of $D_3$.
\begin{proposition} \label{d3}
  Every enumeratively relevant component of $D_H(\P^N,I,J)$ in
  $\mathcal{D}_3 (\P^N,I,J)$  is one of

1. $\mathcal{D}_3 (\P^N,I-1,J+1)$ if $\sum_{j=1}^{J}(N-2-a_j) +
   N-2-b_I \leq  3N + 3$, 

2. $\mathcal{D}_3 (\P^{N-1}, I, J \geq 3)$, or

3. $\mathcal{S}(\P^N, 0,1,1,1,I_1, I \backslash I_1,J \geq 9)$.

\noindent Each of these occurs with multiplicity 1.
\end{proposition}

\noindent {\bf Proof:} This is a complete list of enumeratively
relevant components. If the surface represented by a general point of
a component is irreducible or contained in $H$, we already argued that
cases 1 or 2 must hold. If the surface is reducible with at least one
component in $H$ and one component outside, then the component in $H$
must be a plane. (Note that a trivial dimension count excludes the
possibility that any component is a non-reduced plane.) The component
outside is a possibly reducible quadric surface.  Further specializing
a smooth quadric strictly decreases the dimension and there are no new
contributions from the choice of the limit of the hyperplane section
in $H$. Finally, by Lemma \ref{reduction} the surfaces do not
intersect the intersection of two of the linear spaces $\Omega^i$ or
$\Sigma^j$.  We conclude that we have the complete list.

To prove that the components occur with multiplicity 1, we can assume
that $I=1$, $\Omega^I = \P^N$, and $\Sigma^j = H$.  Using Lemma
\ref{reduction} repeatedly we conclude the proposition in general.
Next by taking a general projection, we can assume that $N=3$. Now the
argument is an easy deformation argument. We will show that as we move
$p_I$ out of the plane there is a first order deformation of the limit
cubic which contains the deformation of the point. This suffices to
conclude that the multiplicity is 1.

Let $(X_0, X_1, X_2, X_3)$ be coordinates on $\P^3$. Suppose $H$ is
defined by $X_0 = 0$. We will write the deformation down for the case
3, as the others are even easier. We can assume that the plane and
quadric pair is given by $X_0 Q(X_0, X_1, X_2 ,X_3)$. Suppose the limit
point is $p= (0,0,0,1)$. We take the deformation $p_{\epsilon} =
(\epsilon, 0, 0,1)$ of the point away from $X_0$. There are already 9
points in $H$. Those define a unique cubic on $X_0$ given by
$C(X_1,X_2,X_3)$. We can assume that $p$ does not lie on $C$. We need
to write a first order deformation of the surface which vanishes on
$C$ and contains $p_{\epsilon}$ to first order:
$$X_0 Q(X_0, X_1, X_2, X_3) - \epsilon \left[ \frac{Q(0,0,0,1)}{C(0,0,1)}
C(X_1,X_2,X_3) + X_0 Q' \right] $$
where $Q'$ is any quadric. $\Box$  \smallskip

Proposition \ref{d3} reduces the problem of counting $D_3$ to counting
plane and quadric pairs meeting in a conic or counting $D_3$ in one
lower dimensional projective space with incidence conditions on a
hyperplane section $E$. 

The techniques in \cite{coskun1:degenerations} readily provide a
solution of the first problem. The incidence conditions on $E$ are
simply incidence conditions on the plane. By Schubert calculus we can
express these conditions in terms of multiples of Schubert cycles. The
problem reduces to counting quadric surfaces whose span contains a
plane satisfying Schubert conditions.

The second problem requires us to describe the enumeratively relevant
components of $D_{\Pi}(\P^N, 3, I, J , J')$ in $D_3(\P^N, 3, I, J ,
J')$, $J' < 4$.

1. The components where $E$ and the surface remain outside $\Pi$ are
of the form $D_3 (\P^N, 3, I, J-1, J'+1)$ if $J' <2$ and $D_3 (\P^N,
3, I, J-1, J'+1, \O(1))$ if $J'=2$. The $\O(1)$ condition on an
elliptic cubic simply means that the intersection points of the curve
with the linear spaces are collinear.

2. If $\Lambda^J$ has dimension $N-2$ and one of the $\Gamma^{j'}$
also has dimension $N-2$, then $E$ can meet $\Lambda^J \cap
\Gamma^{j'}$. If the linear spaces do not have codimension two in
$\P^N$, by Lemma \ref{reduction} the locus where $E$ meets their
intersection is not a divisor.

3. If $E$ lies in $\Pi$ and the surface is irreducible and not in
$\Pi$, then we get the component $D_3 (\P^N, I, J + J')$ when $J'=2$
or $D_3 (\P^N, I, J + J', \O(1))$ when $J' =3$ where the linear spaces
meeting $E$ now are the intersections of the linear spaces with $\Pi$.

4. If the curve breaks, but the surface does not break, then there is
a line in $\Pi$ and a conic not in $\Pi$ meeting it twice.  $D_3(\P^N,
1 +2 =3, I,J,J')$ is the divisor corresponding to this situation.
This is an enumeratively relevant divisor when $J'=3$. The locus where
the conic is more special is a sublocus of this one and hence does not
form a divisor.

There are two other possibilities: the surface can lie in $\Pi$ or
both the surface and $E$ can break.  Neither of these loci give
divisors. The first locus is a specialization of 3. A simple dimension
count excludes the second locus. \smallskip

\noindent {\bf Claim:} All these components occur with multiplicity 1
in $D_{\Pi}(\P^N, 3, I, J , J')$ \smallskip

\noindent {\bf Proof:} By repeatedly applying Lemma \ref{reduction},
we can assume $I=0$.  There is a smooth morphism from
$D_3(\P^N,3,0,J,J')$ to $\overline{M}_{1,J+J'} (\P^N, 3)$ given by
sending the surface curve pair to the stable map which embeds the
curve into $\P^N$. This morphism extends to a general point of the
divisors listed above. We note that in case 2, the map has a
contracted rational component. Both of the spaces are smooth of the
expected dimension at the general points of the listed components and
at their images. To show that the morphism is smooth, it suffices to
check that the Zariski tangent space to the fibers have the expected
dimension. The Zariski tangent space is given by $H^0 (D_3,
N_{D_3/\P^N}(-E))$. By Lemma \ref{counting3} since the latter bundle
does not have any $h^1$, it follows that the morphism is smooth. The
claim is a consequence of Theorem 6.3 \cite{vakil:rationalelliptic}
since $D_{\Pi}(\P^N,3, I,J,J')$ here is the pull-back of $D_H$ in the
notation there by the smooth morphism. $\Box$ \smallskip

Finally, observe that if the data $\tilde{I},\tilde{J}$ differs from
$I,J$ by either including another $N-2$ dimensional linear space to
$I$ or by an $N-1$ dimensional linear space to $J$, then the number of
the new surfaces is the degree of the surface in the case of $I$ and
the degree of the curve in the case of $J$ times the number because of
the choice of the point. See Proposition 6.2
\cite{vakil:rationalelliptic}.

Counting $D_3$ containing a conic or a line is similar, but easier.
\cite{vakil:rationalelliptic} (see \S 7.7-7.9) demonstrates how to
count curves with the $\O(1)$ condition.  Consequently, we can always
count $D_3$ using the degeneration method. As a corollary, we can
count pairs $D_3$ with an elliptic cubic curve which satisfies
incidence conditions with linear spaces and a divisorial condition.
Counting $D_3$ incident to linear spaces can be solved by classical
means, however, at each stage this requires working out the cohomology
ring of a new parameter space. In addition, it is considerably more
difficult to count $D_3$ with curve conditions which satisfy
incidences and divisorial conditions by classical methods. \smallskip

\noindent {\bf Counting $D_4$.} For simplicity we  discuss how to
count $D_4$ in $\P^N$ when at least 4 of the linear spaces are points,
one has dimension $k \leq N-4$ and at least $5+k$ have dimension $N-4$ or
less. We first study the case $N=4$.  The case $N>4$ easily reduces to
it.

Specialize $6$ of the points to a hyperplane $H$ keeping 4 of the
remaining points outside $H$ at all times. Specialize the rest of the
conditions to $H$ in order of increasing dimension. The enumeratively
relevant components of $D_H(\P^4,I,J)$ are one of

1. $D_4 (\P^4, I-1,J+1)$ or 

2. $S(\P^N,1,1,1,1,I_1, I \backslash I_1, J \geq 8)$. 

\noindent These occur with multiplicity 1. 

Since the surface cannot lie in $H$, the first divisor covers all the
possibilities where the surface is irreducible. By Proposition
\ref{junk} and a simple inspection $D_H$ does not contain any
components whose general member parametrizes surfaces with more than 2
components. The only component of $D_H$ where the components of the
limit surface have degree 2 is the divisor 2.  If the components have
degrees 1 and 3, then the cubic spans $\P^4$. If the cubic spanned
$\P^3$, it could meet at most 6 of the points, but the remaining 4 do
not lie on a plane. By Lemmas \ref{conicH} and \ref{glue} there cannot
be a plane and cubic scroll pair either.

To see that they occur with multiplicity 1 we argue as in the case
of $D_3$ by constructing a suitable first order deformation. Since
case 1 is easier, we just write the explicit deformation in case 2.
Without loss of generality we can take the the quadrics given by $X_0
X_1, X_0^2 + X_1^2 + X_2^2 + X_3^2 + X_4^2$. Let $H$ be $X_0 = 0$. The
limit hyperplane section is an elliptic quartic $E$, so it is the
intersection in $X_0=0$ of $X_1^2 + X_2^2 + X_3^2 + X_4^2$ with a
quadric $Q$. We can assume $I=1$ and that the limit point is
$(0,i,0,0,1)$. We deform the point away from $X_0$ by considering
$(\epsilon, i , 0,0,1)$. We need a deformation which vanishes on $E$
and contains the deformation of the point to first order:
\begin{equation*}
\begin{split}
  & X_0 X_1 - \epsilon \left[\frac{i}{Q(i,0,0,1)} Q(X_1, X_2, X_3, X_4) + X_0
  L_1(X_0, X_1, X_2, X_3, X_4) \right]
  \\
  & X_0^2 + X_1^2 + X_2^2 + X_3^2 + X_4^2 + \epsilon X_0 L_2 (X_0,
  X_1, X_2, X_3, X_4)
\end{split}
\end{equation*}  
where $L_i$ are linear forms.

We thus, reduce the problem to counting pairs of quadric surfaces
meeting in a conic where the one in $H$ contains an elliptic quartic
$E$ satisfying incidences. This problem requires us to describe the
enumeratively relevant components of $D_{\Pi}(\P^3;1,1;I,J,J')$ in
$\mathcal{S}(\P^3;1,1;I,J,J')$. We specialize the conditions on $E$
onto a plane by bringing 3 of the 6 points onto the plane keeping 3 of
the points outside and then bringing the rest in order of increasing
dimension. There are at least 8 linear spaces meeting $E$ and at least
6 of them are points. When we specialize $\Gamma^J$

1. $E$ can meet $\Pi$ along $\Gamma^J$. We get
$\mathcal{S}(\P^3;1,1;I,J-1,J'+1)$.

2. If one of the linear spaces in $\Pi$ is already a line and the
condition we specialize is a line, then $E$ can meet their
intersection point. This gives rise to
$\mathcal{S}(\P^3;1,1;I,J-1,J')$.

3. $E$ can break into a conic in $\Pi$ and a conic outside meeting
it twice. This gives rise to $S_{1,1} (\P^3;2+2 = 4, I,J,J')$. In
this case the enumerative problem becomes trivial because the
conditions have to split between $\Pi$ and the plane spanned by the
points outside $\Pi$. We are reduced to counting quadrics subject to
point conditions.

4. The curve can break into a line in $\Pi$ and a twisted cubic
meeting it twice. This gives rise to $S_{1,1} (\P^3;1+3 = 4,
I,J,J')$. There are a finite number of these since the curve has at
least 2 more conditions other than the 6 points
(\cite{vakil:rationalelliptic}). The curve imposes 8 linear conditions
on quadric surfaces, so this case is also very easy to count.

This is a complete list because an easy dimension count shows that
there are no components of $D_{\Pi}$ where the curve has more than two
components and the point conditions imposed on the quadric preclude
the possibility of its breaking into two planes. In particular, the
limit of $E$ cannot be the union of a line with an elliptic cubic. \smallskip

\noindent {\bf Claim:} The divisors above occur with multiplicity
1. \smallskip

\noindent {\bf Proof:}  There is a rational morphism from
$S_{1,1}(\P^3, 4, I,J,J')$ to $\overline{M}_{1,J+J'} (\P^3,4)$ by
sending the marked elliptic curve to the stable map which embeds it in
$\P^3$. This map is well defined on an open set of each of the
divisors and it gives a smooth morphism. The proof is identical to the
case of $D_3$ except instead of using Lemma \ref{counting3}, we use
Lemma \ref{counting4}. $\Box$ \smallskip

Finally, to reduce the more general case to the case $N=4$ specialize
the conditions to a hyperplane containing the $4$ points. Suppose by
induction that we can count $D_4$ in $\P^{N-1}$ satisfying the
analogous conditions and all the degenerations of that case. If the
surface lies in the hyperplane at any stage, then we are reduced to a
subcase of the case $N-1$. If the surface breaks, the argument for
$N=4$ case shows that it must be two quadrics and the same analysis
applies to the quadric in $H$. \smallskip

\noindent {\bf Remark:} For concreteness let us assume $N=4$. Almost
the same argument holds if we required only 8 of the linear spaces to
be points.  In this case the surface can break into a cubic scroll
union a plane in $H$, but these numbers can be determined using
\cite{coskun1:degenerations}.  When exactly 8 of the linear spaces are
points, there are cases when we have to count cubic scrolls meeting a
line twice.  Although this problem is not explicitly addressed in
\cite{coskun1:degenerations}, a simple modification of the argument
for counting twisted cubics meeting a line twice
(\cite{vakil:rationalelliptic} \S 7.5) works.

If only 7 of the linear spaces are points, provided we specialize all
the point but one to $H$ first, the only new limit is $D_3$ union a
plane. These can be counted but requires studying the degenerations of
$D_3$ with an elliptic quartic on it. This case is marginally
harder than counting pairs $D_3$ and an elliptic cubic.

If we relax the conditions a little further, we open Pandora's box.
Singular surfaces appear in the limit and they can have tangencies
with $H$. Identifying the divisors and their multiplicities becomes
difficult. Unlike the simple examples here the multiplicities are not
in general 1. \smallskip

\section{Gromov-Witten Invariants of $X^N$}

In this section we study the relation between the Gromov-Witten
invariants of $X^N$ and the enumerative geometry of $D_n$. We prove
that when $n=3$ and in many cases when $n=4$, the Gromov-Witten
invariants involving incidence to linear spaces are enumerative.
However, when $n \geq 5$, most are not enumerative.  \smallskip

\noindent {\bf Gromov-Witten Invariants.} The Kontsevich spaces of
stable maps possess a virtual fundamental class
$[{\overline{M}_{0,m}(X,\beta)}]^{\text{virt}}$ of the expected
dimension $$
\dim X- K_X \cdot \beta + m -3. $$
They are equipped with
$m$ evaluation morphisms $\rho_1, \cdots, \rho_m$ to $X$, where the $i$-th
evaluation morphism takes the point $[C, p_1, \cdots, p_m, \mu]$ to
the point $\mu(p_i)$ of $X$. Given classes $\gamma_1, \cdots,
\gamma_m$ in the Chow ring $A^* X$ of $X$, we obtain a class

\begin{displaymath}
\rho_1^*(\gamma_1) \cup \cdots \cup \rho_m^* (\gamma_m)
\end{displaymath}
in $\overline{M}_{0,m}(X, \beta)$.  We can evaluate its homogeneous
component of top dimension on
$[\overline{M}_{0,m}(X,\beta)]^{\text{virt}}$ to obtain a number
$I_{\beta}(\gamma_1, \cdots, \gamma_m)$ 
called the {\it Gromov-Witten invariant}. Explicitly,

\begin{displaymath}
I_{\beta}(\gamma_1, \cdots, \gamma_m) =
\int_{[\overline{M}_{0,m}(X,\beta)]^{\text{virt}}} \rho_1^*(\gamma_1)
\cup \cdots \cup   \rho_m^*(\gamma_m).  
\end{displaymath} 
{\bf Notation.} Let $\Gamma_{\Lambda^a} \subset X^N$ denote the
variety of conics in $\P^N$ incident to a linear space $\Lambda^a$ of
dimension $a$. Let $\gamma_a$ denote the cohomology class of
$\Gamma_{\Lambda^a}$.

\begin{definition}
  We call a Gromov-Witten invariant $I_{d_n} (\gamma_{a_1}, \cdots,
  \gamma_{a_m})$ of $X^N$ {\bf enumerative} if for a general set of
  linear spaces $\Lambda^{a_i}$ the only stable maps $(C,p_1, \cdots,
  p_m; \mu)$ in $\overline{M}_{0,m}(X^N,d_n)$ with $\mu_*[C] = d_n$
  and $\mu(p_i) \in \Gamma_{\lambda^{a_i}}$ are injective maps from
  irreducible source curves whose images coincide with a curve of
  conics on a smooth $D_n$. We call these maps {\bf enumerative maps}.
\end{definition}
{\bf Remark.} $I_{d_n} (\gamma_{a_1}, \cdots, \gamma_{a_m})$ is
non-zero only when $$\sum_{i=1}^m (N-1-a_i) = N(n+1) - n + 10 +
m.$$ We always assume that this equality holds and that $a_i < N-2$.

We say a stable map $\mu$ to $X^N$ {\it sweeps out a variety} $V
\subset \P^N$ if the set-theoretic image of the projection of the
universal conic over the image of $\mu$ to $\P^N$ is $V$. If $\mu$
restricted to an irreducible component $C_i$ of the domain curve $C$
sweeps out an irreducible variety $V \subset \P^N$, we say that $(C_i,
\mu_{|C_i})$ {\it sweeps out $V \ $ k times} if the projection from
the universal conic is a generically finite morphism of degree $k$.

If the obstructions for a stable map vanish, then the virtual
fundamental class coincides with the usual one. Lemma 1.1 in
\cite{gathmann:blowup}, which we reproduce for the reader's
convenience, states a local version.

\begin{lemma}\label{Gathmann}
  If $h^1(C, \mu^* T_X)= 0$ for $(C,p_1, \cdots, p_m; \mu) \in
  \overline{M}_{0,m}(X, \beta)$, then $(C,p_1, \cdots, p_m; \mu)$ lies
  in a unique component $Z$ of $\overline{M}_{0,m}(X, \beta)$ of
  dimension equal to the virtual dimension. Moreover,
\begin{displaymath}
[\overline{M}_{0,m}(X, \beta)]^{\mathrm{virt}} = [Z] + R
\end{displaymath}
where $R$
is a cycle supported on the union of the components other than $Z$.
\end{lemma}
Since the normal bundle of a curve $C$ in $X^N$ corresponding to a fixed
conic class on a smooth $D_n$ is generated by global sections, the
standard exact sequence
$$0 \rightarrow T_C \rightarrow \mu^* T_{X^N} \rightarrow N_{\mu}
\rightarrow 0$$
implies that for an enumerative map $h^1(C, \mu^*
T_X)= 0$. By Lemma \ref{Gathmann} an enumerative Gromov-Witten
invariant of $X^N$ is equal to the ordinary scheme-theoretic
intersection of the cycles $\rho_i^* \Gamma_{\Lambda^{a_i}}$ on the
component of enumerative maps, where $\Lambda^{a_i}$ are general
linear spaces. An enumerative map has no automorphisms, so the cycles
intersect at smooth points of this component. Since by Kleiman's
Transversality Theorem they intersect transversely, we conclude

\begin{proposition}\label{gumba}
  Let $\# D_n(\gamma_{a_1}, \cdots, \gamma_{a_m})$ be the number of
  Del Pezzo surfaces incident to general linear spaces of dimension
  $a_i$. Let $R_n$ denote the number of distinct conic classes on
  $D_n$.  Then, for an enumerative Gromov-Witten invariant of $X^N$ we
  have the equality
  $$I_{d_n} (\gamma_{a_1}, \cdots, \gamma_{a_m}) =  R_n \cdot
  \# D_n(\gamma_{a_1}, \cdots, \gamma_{a_m})$$
\end{proposition}   
\noindent {\bf Remark.} Lemma \ref{conics} implies that  $R_3=27,
R_4= 10, R_5= 5, R_{6}= 3$, $R_7= 2$ and $R_8 = 1$.
\smallskip

\noindent {\bf Non-enumerative $I_{d_n}$} The Gromov-Witten invariants
$I_{d_n} (\gamma_{a_1}, \cdots, \gamma_{a_m})$ are in general not
enumerative. The following proposition constructs non-enumerative
examples for $5 \leq n \leq 7$. When $n=8$, the invariants are not
only non-enumerative, but the conics that sweep a $D_8$ and the conics
that sweep the anti-canonical embedding of $\P^1 \times \P^1$ both
contribute to them. Hence, the invariants are not well-suited for
enumerative calculations.

\begin{proposition}\label{nonenumerative}
  The Gromov-Witten invariant $I_{d_n} (\gamma_{a_1}, \cdots,
  \gamma_{a_m})$ of $X^n$ is not enumerative when
\begin{enumerate}
\item  $n=5$ and $\sum_{i=1}^m (2-a_i) \leq 16$,
  
\item $n=6$ and there exists a partition of the numbers $a_i$ into
  $(b_s)_{s=1}^S$ and $(c_t)_{t=1}^T$ such that
   $$
  \sum_{s=1}^S (3-b_s) \leq 12, \ \sum_{t=1}^T (2- c_t) \leq 10, \ 
  \sum_{s=1}^S (4-b_s) = 36, \ \sum_{t=1}^T (4-c_t) = 10, $$

\item  $n=7$ and $\sum_{i=1}^m (3-a_i) \leq 30$. If equality holds we also
  assume that at least one $a_i <3$ is not equal to $1$.
\end{enumerate}
\end{proposition}

\begin{figure}[htbp]
\begin{center}
\epsfig{figure=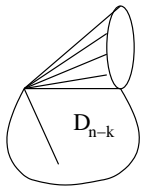}
\end{center}
\caption{The surface $R(k,n)$.}
\label{rnk}
\end{figure}

\noindent {\bf Proof:} For each case we need to construct a
stable map to $X^N$ that satisfies the incidences, but does not sweep
out a smooth $D_n$. We already encountered the additional components
of $\overline{M}_{0,m}(X^N,d_n)$ in the proof of Proposition
\ref{junk}. Take a smooth Del Pezzo surface $D_{n-k}$, $n-k \geq 3$,
and a rational cone of degree $k$ containing a line $l$ of $D_{n-k}$
and meeting a different line $l'$ of $D_{n-k}$ incident to $l$ at its
vertex. Denote this surface by $R(k,n)$ (see Figure \ref{rnk}). The
connected curve of conics corresponding to the conics in the class
$[l]+[l']$ on $D_{n-k}$ and the union of $l'$ with the lines on the
cone has class $d_n$ in $X^N$.
\noindent The 
dimension lemmas in \S 5 imply that the dimension of $R(k,n)$ is
$N(n+1) - n + 10 + k-2$. When $k>1$, this is at least the dimension of
$D_n$. Under the hypotheses of the proposition it is easy to see that
there are $R(2,5)$, $R(2,6)$ and $R(3,7)$ meeting a general set of
$\Lambda^{a_i}$ in $\P^5, \P^6$ and $\P^7$ in the three cases,
respectively. The same construction also provides non-enumerative
examples of $I_{d_n}$ when $N > n$. $\Box$ \smallskip

\begin{theorem}\label{enumerative}
  The Gromov-Witten invariant $I_{d_n} (\gamma_{a_1}, \cdots,
  \gamma_{a_m})$ of $X^N$ is enumerative when 
\begin{enumerate}
\item $n=3$ or 
\item $n=4$ and $\ \sum_{i=1}^m (N-3-a_i) > 4 (N-3)$ .
\end{enumerate}
\end{theorem}

\noindent {\bf Proof:} Let $(C,p_1, \cdots, p_m;\mu)$ be a stable map
to $X^N$ in the class $d_n$ such that $\mu(p_i) \in
\Gamma_{\Lambda^{a_i}}$ for general linear spaces $\Lambda^{a_i}$. Let
$V$ be the variety swept out by $\mu$. Then $V$ must meet the linear
spaces $\Lambda^{a_i}$. The strategy of the proof is to use Theorem
\ref{classification} to prove that $V$ is a smooth $D_n$ under our
hypotheses.  \smallskip

\noindent {\bf Claim 1:} The variety swept out by a connected curve in
the class $d_n$ spans at most $\P^n$. \smallskip

If the curve in $X^N$ is connected, then the threefold in $\P^N$ swept
by the planes of the conics is connected in codimension 1. Since its
degree is bounded by $n-2$, its span (which contains $V$) can be at
most $\P^n$.  \smallskip

\noindent {\bf Claim 2:} $V$ is a variety of pure dimension 2. \smallskip

If $C$ is an irreducible curve in $X^N$, the variety it sweeps in
$\P^N$ does not have to be of pure dimension 2. It can be of pure
dimension 1 or it can have a component of dimension 1 and a component
of dimension 2. 

If a non-constant family in $X^N$ sweeps a curve, then the curve is a
line $l$ and the conics are non-reduced conics whose set-theoretic
support is $l$. Such a component sweeps a surface of degree 0. Since
the degree of the surface swept by the conics in the class $d_n$ is
non-zero, $V$ must contain a surface component. Since the image of
$\mu$ is connected if it contains a component of the type just
described, the line $l$ must lie on a surface component.

\begin{figure}[htbp]
\begin{center}
\epsfig{figure=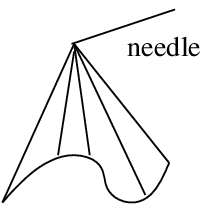}
\end{center}
\caption{Needles.}
\label{needle}
\end{figure}

If an irreducible curve of conics sweeps out a variety which has a
dimension 1 and a dimension 2 component, then the conics consist of
the lines on a cone union a line $l$ meeting the cone at its vertex.
We refer to lines like $l$ as {\it needles.} (See Figure \ref{needle}).

 If $V$ does not have pure dimension 2, then it must contain
a needle. Since the image of $\mu$ is connected, if the needle is not
contained in a surface component, then every component of the image of
$\mu$ must have the same line as a needle.  A stable map onto such a
curve cannot be in the class $d_n$. The class of a stable map which
sweeps a rational cone of degree $r$ with a needle $k$ times is $-kr\ 
a + kr \ b$. Hence, $V$ has pure dimension 2.  \smallskip

\noindent {\bf Claim 3:} If $n=3$ (resp. $4$), $V$ is an irreducible
surface of degree $3$ (resp. $4$). \smallskip

Since meeting linear spaces in general position impose independent
conditions, we can check this claim by a naive dimension count.
$V$ has degree at most $n$ and spans at most $\P^n$. 

When $n=3$, the linear spaces impose $19 + 4(N-3)$ conditions on $V$.
If $V$ has degree less than 3, then it is either a plane or a quadric
surface, hence $V$ has dimension at most $9 + 4(N-3)$.  We conclude
that $V$ has degree 3. This concludes the proof when $n=3$ because any
surface of degree $3$ in $\P^3$ is a specialization of smooth $D_3$
surfaces. The variety of singular cubic surfaces has codimension 1 in
the space of cubic surfaces. If the linear spaces are general, the
only cubic surfaces that satisfy all the incidences will be smooth
surfaces. We conclude that the image of $\mu$ must contain a curve of
conics on $D_3$, but then the image must coincide with it. Also
observe that $\mu$ cannot have any contracted components, since
otherwise $V$ would have to meet $\Lambda^{a_i}$ and $\Lambda^{a_j}$
for some $i,j$ along the same conic. Dimension considerations exclude
this possibility.

When $n=4$, there are more cases to consider. It suffices to carry out
the dimension counts when $N=4$. Since the dimension of cubic surfaces
in $\P^4$ is bounded by 23, the degree of $V$ must be 4.

A surface swept by a connected curve in $X^N$ is connected in
codimension 1 by lines or conics except when it contains  cones
meeting only along their vertices.

The dimension of four-tuples of planes in $\P^4$ or triples of a
quadric surface and two planes in $\P^4$ is bounded by 25. Hence, if
$V$ is reducible, then it is either the union of two quadrics or the
union of a cubic surface and a plane.  If the quadrics share a common
line or conic or if one of them is a cone, then their dimension is
strictly less than 26. The choice of pairs of a plane and a cubic that
spans $\P^4$ is bounded by $24$.  If the cubic spans only $\P^3$ and
is ruled by lines, it is either the projection of a cubic scroll or a
cone over an elliptic curve. In either case the dimension of the
choice of pairs of such a cubic and a plane is bounded by 23 (\S 5.1).
If the cubic is not ruled by lines, then it contains only a one
parameter family of conics and the plane must meet it along a line or
a conic.  By \S 5.1 the dimension of these pairs is also strictly less
than 26. We conclude that $V$ is irreducible.

Theorem \ref{classification} classifies irreducible quartic surfaces
that span $\P^4$. The projection of a scroll of degree 4 in $\P^5$ and
of a Veronese surface have dimensions bounded by $23$ and $16$,
respectively. The other possibilities are degenerations of $D_4$,
hence have dimension strictly smaller than $26$. Finally by the
assumption that $\sum_{i=1}^m (N-3-a_i) > 4 (N-3)$, no $\P^3$ meets
all the linear spaces $\Lambda^{a_i}$. We conclude that $V$ must span
$\P^4$. This completes the proof that $V$ is a smooth $D_4$. The rest
of the argument is identical to the previous case.  $\Box$ \smallskip

\noindent {\bf Remark:} When $n=4$, can we remove the assumption on
$a_i$? Any counterexample must arise from an irreducible quartic
surface in $\P^3$. Suppose we take a quartic surface in $\P^3$ with a
double line. The dimension of such surfaces is $25$. The conics that
are residual to the double line $l$ in the pencil of planes containing
it give us a curve $C$ in $X^3$.  The class of $C$ is $2a + b$, not
$d_4$. Suppose now we choose a more special quartic surface $S$ so
that $C$ contains a non-reduced conic whose set theoretic support is a
line $m$. Take the curve $C'$ of non-reduced conics whose set
theoretic supports are $m$, but whose planes rotate once about $m$ in
$\P^3$.  The union of $C$ and $C'$ now are in the class $d_4$. The
dimension of surfaces $S$ in $\P^3$ with the required property is
$22$. We conclude that Theorem \ref{enumerative} is the sharpest we
can hope for when $n=4$. \smallskip

The table below gives examples of Gromov-Witten invariants of $X^N$ we
can calculate using Theorem \ref{enumerative}, Proposition \ref{gumba}
and the degeneration method in \S 5.  We use the short-hand
$I_{d_n}(a_1^{r_1}, \cdots, a_k^{r_k})$ to denote the Gromov-Witten
invariant of $X^N$ in the class $d_n$ incident to $r_i$ cycles of
conics meeting linear spaces of dimension $a_i$. \bigskip

\begin{tabular}{|cl|cl|} \hline
$N=3$ & $I_{d_3}(0^{19})=27$ \ & $N=4$ & $I_{d_4}(0^{13})= 10$ \\ \hline
$N=4$ & $I_{d_3}(0^4,1^{15})=27$ \ & $''$ &  $I_{d_4}(0^{12},1^2)= 40$ \\
\hline
$''$ & $I_{d_3}(0^3,1^{17})=972$ \ & $''$ &  $I_{d_4}(0^{11},1^4)= 320$ \\
\hline
$''$ & $I_{d_3}(0^2,1^{19})= 21303$ \ & $''$ &  $I_{d_4}(0^{10},1^6)= 3200$ \\
\hline
$N=5$ & $I_{d_3}(0^2,1^{4}, 2^{13})= 54$ \ & $''$ & $I_{d_4}(0^{9},1^8)=
33280$ \\ \hline
$''$ & $I_{d_3}(0^2,1^3,2^{15}) = 1863$ \ & N=5 & $I_{d_4}(0^{4},1^9,2)=
240$ \\ \hline
\end{tabular}

\bibliographystyle{math} 
\bibliography{math} 

\bigskip

\noindent Mathematics Department, Harvard University, Cambridge, MA
02138 

\noindent E-mail: coskun@math.harvard.edu \end{document}